\documentclass[12pt]{amsart}
\usepackage{amssymb}
\usepackage{epsf}

\newtheorem{thm}{Theorem}[section]
\newtheorem{prop}[thm]{Proposition}
\newtheorem{lem}[thm]{Lemma}
\newtheorem{cor}[thm]{Corollary}
\newtheorem{definition}[thm]{Definition}
\newenvironment{defn}{\begin{definition}\sl}{\end{definition}}
\newtheorem{remark}[thm]{Remark}
\newenvironment{rem}{\begin{remark}\rm}{\end{remark}}
\newtheorem{example}[thm]{Example}
\newenvironment{ex}{\begin{example}\rm}{\end{example}}

\newcommand{\QED}{
\setlength{\unitlength}{1.0pt}%
\begin{picture}(7.5,7.5)
\put(0,-5){\rule{2.5pt}{2.5pt}}
\put(0,-2.5){\rule{2.5pt}{2.5pt}}
\put(0,0){\rule{5pt}{2.5pt}}
\put(0,2.5){\rule{7.5pt}{2.5pt}}
\end{picture}\vspace{10pt}}

\begin{document}

\title[Skew Schubert functions]{Skew Schubert functions and the\\ Pieri
formula for flag manifolds} 

\author{Nantel Bergeron \and Frank Sottile}

\address{Department of Mathematics and Statistics\\
        York University\\
        North York, Ontario M3J 1P3\\
	CANADA}
\email[Nantel Bergeron]{bergeron@mathstat.yorku.ca}
\urladdr[Nantel Bergeron]{http://www.math.yorku.ca/bergeron}
\address{Department of Mathematics\\
        University of Toronto\\
        100 St.~George Street\\
	Toronto, Ontario  M5S 3G3\\
	CANADA}
\email[Frank Sottile]{sottile@math.toronto.edu}
\urladdr[Frank Sottile]{http://www.math.toronto.edu/\~{}sottile}
\date{\today}
\thanks{First author supported in part by NSERC and CRM grants}
\thanks{Second author supported in part by NSERC grant  OGP0170279,  
NSF grant DMS-9022140, and CRM}

\begin{abstract}
We show the equivalence of the Pieri formula for flag manifolds and
certain identities among the structure constants,
giving new proofs of both the Pieri formula and of these identities.
A key step is the association of a symmetric function to a finite 
poset with labeled Hasse diagram satisfying a symmetry condition.
This gives a unified definition of skew Schur functions, Stanley 
symmetric function, and skew Schubert functions (defined here).
We also use algebraic geometry to show the coefficient of a monomial in a
Schubert polynomial counts certain chains in the Bruhat order,
obtaining a new combinatorial construction of Schubert polynomials.
\end{abstract}

\maketitle

\section*{Introduction}

A fundamental open problem in the theory of Schubert
polynomials is to find an analog of the Littlewood-Richardson
rule.
By this, we mean a bijective description of the structure
constants for the ring of polynomials with respect to its basis of Schubert
polynomials.
Such a rule would express the intersection form in the cohomology
of a flag manifold in terms of its basis of Schubert classes.
Other than the classical Littlewood-Richardson rule (when the Schubert
polynomials are Schur symmetric polynomials) little is known. 

Using geometry, Monk~\cite{Monk} established a formula for multiplication by
linear Schubert polynomials (divisor Schubert classes).
A Pieri-type formula for multiplication by an elementary or 
complete homogeneous symmetric polynomial (special Schubert class) was given
in~\cite{LS82a}, but only recently proven~\cite{Sottile96}
using geometry.
There are now several
proofs~\cite{C-F,Postnikov,Winkel_multiplication,Veigneau},  
some of which~\cite{Postnikov,Winkel_multiplication,Veigneau}
are purely combinatorial.

In the more general setting of multiplication by a Schur symmetric
polynomial, formulas for some structure constants follow from
a family of identities which were proven using geometry~\cite{BS97a}.
Also in ({\em ibid.}) are combinatorial results about intervals in
the Bruhat order which are formally related to these identities.
A combinatorial (but {\em not} a bijective) formula was given
for these coefficients~\cite{BS_monoid} using the Pieri formula,
which gave a direct connection between some of these 
order-theoretic results and identities.

A first goal of this paper is to deduce another
identity~\cite[Theorem~G({\em ii})]{BS97a} 
from the Pieri formula, and also to deduce the Pieri formula
from these identities.
This furnishes a new proof of the Pieri formula, shows its equivalence to
these (seemingly) more general identities,
and, together with the combinatorial proofs of the Pieri formula, 
gives a purely
combinatorial proof of these identities. 

A key step is the definition of a symmetric function 
associated to any finite {\em symmetric labeled poset}, which is a poset 
whose Hasse diagram has edges labeled with integers with a symmetry 
condition satisfied by its maximal chains.
This gives a unified construction of skew Schur functions (for
intervals in Young's lattice of partitions), Stanley symmetric
functions~\cite{Stanley84} (for intervals in the weak order on the
symmetric group), and for intervals in a $k$-Bruhat order,  
{\em skew Schubert functions} (defined in another fashion in \S 1).

In~\cite{LS82b}, Lascoux and Sch\"utzenberger show that if a 
Schubert polynomial is expressed as a univariate polynomial in the first 
variable, then the coefficients are (explicitly determined)
multiplicity-free  sums of Schubert polynomials in the remaining variables.
This may be used to show that Schubert polynomials
are sums of monomials with non-negative coefficients.
We use a cohomological formula~\cite[Theorem 4.5.4]{BS97a} to 
generalize their result, obtaining a similar formula for expressing  a 
Schubert polynomial as a polynomial in {\em any} variable.
This also extends Theorem~C~({\em ii}) of~\cite{BS97a}, which identified
the constant term of this expression.
{}From this, we obtain a construction of Schubert
polynomials purely in terms of chains in the Bruhat order, and a
geometric proof that the monomials which appear in a Schubert
polynomial have non-negative coefficients.
The Pieri formula shows these coefficients are certain intersection numbers,
recovering a result of Kirillov and Maeno~\cite{KM}.

We found these precise formulas in terms of intersection numbers
surprising;
Other combinatorial constructions are either
recursive~\cite[4.17]{Macdonald91} 
and do not give the coefficients, or are expressed in terms of 
combinatorial structures (the  
weak order on the symmetric group~\cite{BJS,FoSt,FK_YB}
or diagrams of permutations~\cite{Kohnert,Be92,Winkel_kohnert_rule})
which are not geometric.
Previously, we believed this non-negativity of monomials had no relation to
geometry. 
Indeed, only monomials of the form $x^\lambda$ with $\lambda$ a partition
are represented by positive cycles, 
other polynomial representatives of Schubert classes~\cite{BGG,De74}
do not have this non-negativity, and polynomial representatives for the other
classical groups~
cannot~\cite{FK_Bn} have such non-negativity. 

This paper is organized as follows.
In Section 1, we give necessary background, define skew Schubert functions, 
and state our main results.
In Section 2, we deduce the Pieri formula from the identities and 
results on the Bruhat order.
In Section 3, we define a symmetric function $S_P$ associated to a 
symmetric labeled poset $P$ and complete the proof of the equivalence 
of the Pieri formula and these identities.
We also show how this construction gives skew Schur and Schubert
functions.
In Section 4, we adapt an argument of Remmel and
Shimozono~\cite{Remmel_Shimozono} to show that, for intervals in the weak
order, this symmetric function is Stanley's symmetric
function~\cite{Stanley84}. 
Finally, in Section 5, we use a geometric result of~\cite{BS97a}
to generalize the result in~\cite{LS82b} and interpret the coefficient of a
monomial in a Schubert polynomial in terms of chains in the Bruhat order.

\section{Preliminaries}

Let ${\mathcal S}_n$ be the symmetric group on $n$ letters and 
${\mathcal S}_\infty:=\bigcup_n{\mathcal S}_n$, the 
group of permutations of ${\mathbb N}$ which fix all but finitely 
many integers.
We let 1 be the identity permutation.
For each $w\in{\mathcal S}_\infty$, 
Lascoux and Sch\"utzenberger~\cite{LS82a} defined a Schubert polynomial
${\mathfrak S}_w \in{\mathbb Z}[x_1,x_2,\ldots]$ with 
$\deg{\mathfrak S}_w = \ell(w)$.
These satisfy the following:
\begin{enumerate}
\item $\{{\mathfrak S}_w\mid w\in{\mathcal S}_\infty\}$ is a 
${\mathbb Z}$-basis for ${\mathbb Z}[x_1,x_2,\ldots]$.
\item 
If $w$ has a unique descent at $k$ ($w(j)>w(j+1)\Rightarrow j=k$), then 
${\mathfrak S}_w = S_\lambda(x_1,\ldots,x_k)$,
where $\lambda_j = w(k+1-j)-k-1+j$.
We write $v(\lambda,k)$ for this permutation
and call $w$ a {\em Grassmannian permutation} with descent $k$.
\end{enumerate}

By the first property, there exist integral structure constants $c^w_{u\,v}$
for  $w,u,v\in{\mathcal S}_\infty$ (non-negative from geometry) defined 
by the identity
\begin{equation}\label{eq:cwuv}
	{\mathfrak S}_u \cdot {\mathfrak S}_v \ =\ 
	\sum_w c^w_{u\,v}\,{\mathfrak S}_w.
\end{equation}
We are concerned with the coefficients $c^w_{u\;v(\lambda,k)}$
which arise when ${\mathfrak S}_v$ in~(\ref{eq:cwuv}) is replaced by the 
Schur polynomial 
$S_\lambda(x_1,\ldots,x_k)={\mathfrak S}_{v(\lambda,k)}$.

It is well-known (see for example~\cite{Sottile96,BS97a}) that 
$c^w_{u\; v(\lambda,k)}\neq0$ only if $u\leq_k w$, where $\leq_k$
is the $k$-Bruhat order (introduced in~\cite{LS83}).
In fact, $u\leq_k w$ if and only if there is some $\lambda$ with
$c^w_{u\; v(\lambda,k)}\neq0$.
This suborder of the Bruhat order has the following characterization:

\begin{definition}[Theorem~A of~\cite{BS97a}]\label{def:1}
{\sl
Let $u,w\in{\mathcal S}_\infty$.
Then $u\leq_k w$ if and only if
\begin{enumerate}
\item $a\leq k<b \Longrightarrow u(a)\leq w(a)$ and $u(b)\geq w(b)$,
\item $a<b, u(a)<u(b)$, and $w(a)>w(b) \Longrightarrow a\leq k<b$.
\end{enumerate}
}\end{definition}

For any infinite subset $P$ of ${\mathbb N}$, the order-preserving bijection 
${\mathbb N}\leftrightarrow P$ and the inclusion 
$P\hookrightarrow {\mathbb N}$ induce a map
$$
\varepsilon_P\ : \ {\mathcal S}_\infty\ \simeq\ {\mathcal S}_P \ 
\hookrightarrow\ {\mathcal S}_\infty.
$$
{\em Shape-equivalence} is the equivalence relation generated by
$\zeta \sim \varepsilon_P(\zeta)$ for $P\subset{\mathbb N}$.

If $u\leq_k w$, let $[u,w]_k$ denote the interval between $u$ and $w$ 
in the $k$-Bruhat order.
These intervals have the following property:
\medskip

\noindent{\bf Order 1} (Theorem~E({\em i}) of~\cite{BS97a}){\bf .}
{\em 
Suppose $u,w,y,z\in{\mathcal S}_\infty$ with $u\leq_k w$,
$y\leq_l z$, and $wu^{-1}$ shape-equivalent to $zy^{-1}$.
Then $[u,w]_k\simeq[y,z]_l$.
Moreover, if $zy^{-1}=\varepsilon_P(wu^{-1})$, then this isomorphism is 
induced by the map $v\mapsto \varepsilon_P(vu^{-1}) y$.}\medskip

This has a companion identity among the structure constants
$c^w_{u\;v(\lambda,k)}$:
\medskip

\noindent{\bf Identity 1} (Theorem~E({\em ii}) of~\cite{BS97a}){\bf .}
{\em 
Suppose $u,w,y,z\in{\mathcal S}_\infty$ with $u\leq_k w$,
$y\leq_k l$, and $wu^{-1}$ shape-equivalent to $zy^{-1}$.
Then, for any partition $\lambda$,
$$
c^w_{u\;v(\lambda,k)}\ =\ 
c^z_{y\;v(\lambda,l)}.
$$
}

This identity was originally proven using geometry~\cite{BS97a}.
In~\cite{BS_monoid}, we showed how to deduce it from Order~1 and the 
Pieri formula for Schubert polynomials.
Here, we use it to deduce the Pieri formula.

By Identity 1, we may define a constant $c^\zeta_\lambda$ 
for any permutation $\zeta\in{\mathcal S}_\infty$ and any 
partition $\lambda$ by $c^\zeta_\lambda=c^w_{u\;v(\lambda,k)}$
for any $u\leq_k w$ with $w=\zeta u$.
We also define the {\em skew Schubert function} $S_\zeta$ by
\begin{equation}\label{eq:skew_Schubert}
S_\zeta\ =\ \sum_\lambda c^\zeta_\lambda S_\lambda,
\end{equation}
where  $S_\lambda$ is the Schur symmetric function~\cite{Macdonald95}.

By  Order~1, we may make the following definition:

\begin{defn}
Let $\eta,\zeta\in{\mathcal S}_\infty$.
Then $\eta\preceq\zeta$ if and only if 
there is a $u\in{\mathcal S}_\infty$ and $k\in{\mathbb N}$ with 
$u\leq_k \eta u\leq_k \zeta u$.
For $\zeta\in{\mathcal S}_\infty$, define
$|\zeta|:= \ell(\zeta u)-\ell(u)$ for any $u,k$ with 
$u\leq_k \zeta u$.  (There always is such a $u$ and $k$, see \S 2.)
\end{defn}

In \S2, $\preceq$ and $|\zeta|$ are given definitions that do not 
refer to $\leq_k$ or $\ell(w)$.

Let $\zeta,\eta\in{\mathcal S}_\infty$.
If we have $\eta\cdot\zeta=\zeta\cdot\eta$ with 
$|\zeta\cdot\eta|=|\zeta|+|\eta|$,
and neither of $\zeta$ or $\eta$ is the identity, 
then we say that $\zeta\cdot\eta$ is the 
{\em disjoint product} of $\zeta$ and $\eta$.
If a permutation cannot be written in this way, then 
it is {\em irreducible}.
It is a consequence of~\cite[\S 3]{BS97a} that
a permutation $\zeta$ factors uniquely into irreducibles
as follows:
Let $\Pi$ be the finest non-crossing partition~\cite{Kreweras} which is
refined  
by the partition given by the cycles of $\zeta$.
For each non-singleton part $p$ of  $\Pi$, let 
$\zeta_p$ be the product of cycles which partition $p$.
Each $\zeta_p$ is irreducible, and $\zeta$ is the disjoint product 
of the $\zeta_p$'s.
See Remark~\ref{rem:cyclic} for a further discussion.
\medskip

\noindent{\bf Order 2} (Theorem~G({\em i}) of~\cite{BS97a}){\bf .}
{\em 
Suppose $\zeta=\zeta_1\cdots\zeta_t$ is the factorization of 
$\zeta\in{\mathcal S}_\infty$ into irreducibles.
Then the map
$(\eta_1,\ldots,\eta_t)\mapsto \eta_1\cdots\eta_t$
induces an isomorphism
$$
[1,\zeta_1]_\preceq\times\cdots\times[1,\zeta_t]_\preceq\ 
\stackrel{\sim}{\relbar\joinrel\longrightarrow}\ [1,\zeta]_\preceq.
$$
}

\noindent{\bf Identity 2} (Theorem~G({\em ii}) of~\cite{BS97a}){\bf .}
{\em 
Suppose $\zeta=\zeta_1\cdots\zeta_t$ is the factorization of 
$\zeta\in{\mathcal S}_\infty$ into irreducibles.
Then 
$$
S_\zeta\ =\ S_{\zeta_1}\cdots S_{\zeta_t}.
$$
}

Theorem G({\em ii}) in~\cite{BS97a} states that if $\zeta\cdot\eta$ is
a disjoint product, then, for all partitions $\lambda$, 
$$
c^{\zeta\cdot\eta}_\lambda\ =\ 
\sum_{\mu,\,\nu} c^\lambda_{\mu\,\nu} c^\zeta_\mu c^\eta_\nu.
$$
Thus we see that 
\begin{eqnarray*}
 S_\zeta \cdot S_\eta   &=& 
\sum_{\mu,\,\nu}  c^\zeta_\mu c^\eta_\nu S_\mu S_\nu\\
&=&  \sum_{\lambda,\,\mu,\,\nu}
c^\lambda_{\mu\,\nu}c^\zeta_\mu c^\eta_\nu S_\lambda\\
&=&
\sum_\lambda c^{\zeta\cdot\eta}_\lambda S_\lambda\ \:=\ \:
S_{\zeta\cdot\eta}.
\end{eqnarray*}
Iterating this shows the equivalence of Theorem G({\em ii}) of~\cite{BS97a}
and Identity 2.
\bigskip

A {\em labeled poset} $P$  is a finite ranked poset together with an integer
label for each cover.
Its Hasse diagram is thus a directed labeled graph with integer labels.
Write $u\stackrel{\mbox{\scriptsize $b$}}{\longrightarrow}w$ 
for a labeled edge in this Hasse diagram.
In what follows, we consider four classes of labeled posets:
\begin{enumerate}
\item[] 
\hspace{-26pt}{\bf Intervals in a $k$-Bruhat order.}
Labeling a cover $u\lessdot_k w$ in the $k$-Bruhat order with $b$, 
where $wu^{-1}=(a,\,b)$ and $a<b$ gives every interval in 
the $k$-Bruhat order the structure of a labeled poset.
\item[] 
\hspace{-26pt}{\bf  Intervals in the $\preceq$-order.}
Likewise, a cover $\eta\prec\!\!\!\!\!\cdot\ \zeta$ in the $\preceq$-order
gives a transposition $(a,\,b)=\zeta\eta^{-1}$ with 
$a<b$.
Labeling such a cover with $b$ gives every interval in this order
the structure of a labeled poset.
Since $[\eta,\zeta]_\preceq \simeq [1,\zeta\eta^{-1}]_\preceq$, it suffices
to consider intervals of the form $[1,\zeta]_\preceq$.
\item[] 
\hspace{-26pt}{\bf Intervals in Young's lattice.}
A cover $\mu\subset\!\!\!\!\!\cdot\ \lambda$ in Young's lattice of partitions 
gives a unique index $i$ with $\mu_i+1=\lambda_i$.
Labeling such a cover with $\lambda_i-i$ 
gives every interval in Young's lattice the structure of 
a labeled poset.
\item[] 
\hspace{-26pt}{\bf Intervals in the weak order.}
Finally, labeling a cover $u\lessdot_{\mbox{\scriptsize\rm weak}}w$ in the
weak order on ${\mathcal S}_\infty$ with the index $i$ of the transposition
$wu^{-1}=(i,\,i{+}1)$ gives every interval in the weak order the structure of
a labeled poset. 
Since, for $u\leq_{\mbox{\scriptsize\rm weak}}w$, 
$[u,w]_{\mbox{\scriptsize\rm weak}}\simeq
[1,wu^{-1}]_{\mbox{\scriptsize\rm weak}}$, it suffices to consider intervals
of the form $[1,w]_{\mbox{\scriptsize\rm weak}}$.
\end{enumerate}

The sequence of edge labels in a (maximal) chain of a labeled poset 
is the {\em word} of that chain.
For a composition $\alpha=(\alpha_1,\ldots,\alpha_k)$ of $m=$ rank$P$,
let $H_\alpha(P)$ be the set of maximal chains in $P$ whose word has 
descent set contained in 
$I(\alpha):=\{\alpha_1,\alpha_1+\alpha_2,\ldots,m-\alpha_k\}$.
We say that $P$ is {\em symmetric} if the cardinality of 
$H_\alpha(P)$ depends only upon the parts of $\alpha$ and not their order.

Each poset in the above classes is symmetric:
For the $k$-Bruhat orders or $\preceq$ order, this is a consequence of
the Pieri formula for Schubert polynomials.
For Young's lattice, this is classical, and for intervals in the 
weak order, it is due to Stanley~\cite{Stanley84}.

We wish to consider skew Young diagrams to be equivalent if they differ by 
a translation.
This leads to the following notion of isomorphism for labeled posets.

\begin{defn}\label{def:lposet}
A map $f:P\rightarrow Q$ between labeled posets is an isomorphism if $f$ is
an isomorphism of posets which preserves the relative order of the edge
labels.
\end{defn}

That is, if $e,e'$ are edges of $P$ with respective labels 
$a\leq a'$, then the edge labels $b,b'$ of $f(e),f(e')$ in $Q$
satisfy $b\leq b'$.
The isomorphisms of Order~1 and Order~2 are isomorphisms
of labeled posets.
We also see that the interval $[\mu,\lambda]_\subset$ in Young's poset 
is isomorphic to the interval $[v(\mu,k),v(\lambda,k)]_k$, 
since the difference between the label of a cover
$v(\alpha,k)\lessdot_k v(\beta,k)$ in the $k$-Bruhat order
and the corresponding cover 
$\alpha\subset\!\!\!\!\cdot\ \beta$ in Young's lattice
is $k+1$.

To every symmetric labeled poset $P$, we 
associate~(Definition~\ref{def:skew}) a
symmetric function $S_P$ which has the following properties:

\begin{thm}\label{thm:skew}
\
\begin{enumerate}
\item If $P\simeq Q$, then $S_P=S_Q$.
\item If $u\leq_k w$, then $S_{[u,w]_k} = S_{wu^{-1}}$, the skew Schubert
function.
\item[2$'$\!.] For $\zeta\in{\mathcal S}_\infty$, 
$S_{[1,\zeta]_\preceq}= S_\zeta$, the skew Schubert function. 
\item Let $\mu\subset \lambda$ be partitions.
Then $S_{[\mu,\lambda]_\subset} = S_{\lambda/\mu}$, the skew Schur function.
\item For $w\in{\mathcal S}_\infty$, we have
$S_{[1,w]_{\rm weak}} = F_w$, the Stanley symmetric function.
\end{enumerate}
\end{thm}

Part 1 is Lemma~\ref{lem:coeff}(2), parts 2, 2$'$, and 3 are proven in 
\S3, and part 4 in \S4.

A labeled poset $P$ is an {\em increasing chain} if it is totally ordered
with increasing edge labels.
A cycle $\zeta\in{\mathcal S}_\infty$ is {\em increasing} if 
$[1,\zeta]_\preceq$ is an increasing chain.
Decreasing chains and cycles are defined similarly.

For any positive integers $m,k$ let $r[m,k]$ denote the 
permutation $v((m,0,\ldots,0),\,k)$ which is the increasing cycle
$(k{+}m,k{+}m{-}1,\ldots,k)$.
It is an easy consequence (see Lemma~\ref{lem:perm_facts}) of the
definitions of $\leq_k$ 
or $\preceq$ that any increasing cycle $\zeta$ of length $m{+}1$
is shape equivalent to $r[m,k]$ and hence $|\zeta|=m$.
Likewise, the permutation $v(1^m,k)$ is the decreasing cycle
$(k{+}1{-}m,\ldots,k,k{+}1)$ and any decreasing cycle of length $m{+}1$ is 
shape equivalent to $v(1^m,k)$ for any $k\geq m$.
Here $1^m$ is the partition of $m$ into $m$ equal parts of size $1$.
Note that 
$$
{\mathfrak S}_{r[m,k]}\ =\ h_m(x_1,\ldots,x_k)
\qquad\mbox{and}\qquad
{\mathfrak S}_{v(1^m,k)}\ =\ e_m(x_1,\ldots,x_k),
$$
the complete homogeneous and elementary symmetric polynomials.

\begin{prop}[Pieri formula for Schubert polynomials and flag
manifolds]\label{Pieri_formula} 
Let $u\leq_k w$ with $m=\ell(w)-\ell(u)$.  
Then
\begin{enumerate}
\item ${\displaystyle c^w_{u\, r[m,k]}\ =\ \left\{
\begin{array}{ll}
1&\ \mbox{if $wu^{-1}$ is the disjoint product of increasing cycles}\\
0&\ \mbox{otherwise.} \end{array}\right.}$
\item ${\displaystyle c^w_{u\, v(1^m,k)}\ =\ \left\{
\begin{array}{ll}
1&\ \mbox{if $wu^{-1}$ is the disjoint product of deceasing cycles}\\
0&\ \mbox{otherwise.} \end{array}\right.}$
\end{enumerate}
\end{prop}

This is the form of the Pieri formula stated in~\cite{LS82a},
as such a disjoint products of increasing (decreasing) cycles are
$k$-{\em soul\`evements droits} (respectively {\em gauches}) for $u$.
By~\cite[Lemma 6]{Sottile96}, $wu^{-1}$ is a disjoint product of increasing
cycles if and only if there is a maximal chain in $[u,w]_k$ with increasing
labels, and such chains are unique.
When this occurs, we write 
$u\stackrel{r[m,k]}{\relbar\joinrel\relbar\joinrel\longrightarrow} w$, where
$m:= \ell(w)-\ell(u)$.
Similarly, $wu^{-1}$ is a disjoint product of decreasing cycles if and only
if there is a maximal chain in $[u,w]_k$ with decreasing labels, 
which is necessarily unique.

Recall that 
\begin{eqnarray*}
H^*(\mbox{\em Flags}({\mathbb C}\,^n))
&\simeq& {\mathbb Z}[x_1,x_2,\ldots]/
\langle {\mathfrak S}_w\mid w\not\in{\mathcal S}_n\rangle\\
&=& {\mathbb Z}[x_1,\ldots,x_n]/
\langle x^{\alpha} \mid \alpha_i\geq n-i,\mbox{ for some } i\rangle.
\end{eqnarray*}
The map defined by ${\mathfrak S}_w \mapsto {\mathfrak S}_{\overline{w}}$,
where $\overline{w}=\omega_0 w\omega_0$, conjugation by the longest element
$\omega_0$ in ${\mathcal S}_n$, 
is an algebra involution on $H^*(\mbox{\em Flags}({\mathbb C}\,^n))$.
If $n\geq k+m$, then this involution shows the equivalence of
the two versions of the Pieri formula.

We state the main results of this paper:

\begin{thm}\label{thm:main_equiv}
Given the results Order 1 and 2 on the $k$-Bruhat
orders/$\preceq$-order, the Pieri formula for Schubert polynomials is
equivalent to the Identities 1 and 2.
\end{thm}

This is proven in \S2 and \S3.

\begin{thm}
If $w\in{\mathcal S}_n$ and $0\leq \alpha_i\leq n-i$ for $1\leq i\leq n-1$,
then the coefficient of
$x_1^{n-1-\alpha_1}x_2^{n-2-\alpha_2}\cdots x_{n-1}^{1-\alpha_{n-1}}$ 
in the Schubert polynomial ${\mathfrak S}_w(x)$ is the number of chains
$$
w \stackrel{r[\alpha_1,1]}{\relbar\joinrel\relbar\joinrel\longrightarrow}
w_1 \stackrel{r[\alpha_2,2]}{\relbar\joinrel\relbar\joinrel\longrightarrow}
\cdots \stackrel{r[\alpha_{n-1},n-1]}{\relbar\joinrel\relbar\joinrel%
\relbar\joinrel\relbar\joinrel\relbar\joinrel\longrightarrow}
\omega_0
$$
between w and $\omega_0$, the longest element in ${\mathcal S}_n$.
\end{thm}

This is a restatement of Corollary~\ref{cor:chain_monomial}.

\section{Proof of the Pieri formula for Schubert polynomials and flag
manifolds} 

Here, we use Identities 1 and 2 to deduce the Pieri formula.
We first establish some combinatorial facts about chains and
increasing/decreasing cycles.

Let $\zeta\in{\mathcal S}_\infty$.  
We give a $u\in{\mathcal S}_\infty$ and $k>0$ such that 
$u\leq_k\zeta u$ and $\zeta u$ is Grassmannian of descent $k$.
Define up$(\zeta):=\{a\mid a<\zeta(a)\}$,
down$(\zeta):=\{b\mid b>\zeta(b)\}$,
fix$(\zeta):=\{c\mid c=\zeta(c)\}$,
and set $k:=\#\mbox{up}(\zeta)$.
If we have
\begin{eqnarray*}
\mbox{up}(\zeta)&=&\{a_1,\ldots,a_k\mid
\zeta(a_1)<\zeta(a_2)<\cdots<\zeta(a_k)\},\\
\mbox{fix}(\zeta)\bigcup\mbox{down}(\zeta)&=&\{b_1,b_2,\ldots\mid
\zeta(b_1)<\zeta(b_2)<\cdots\},
\end{eqnarray*}
and define $u\in{\mathcal S}_\infty$ by
$$
u\ :=\ \left\{\begin{array}{ll}
a_i&\ \mbox{if } i\leq k\\
b_{i-k}&\ \mbox{if } i>k\end{array}\right.,
$$
then $u\leq_k\zeta u$.
Set $w:=\zeta u$.

This construction of $u\in{\mathcal S}_\infty$ 
is Theorem 3.1.5 ({\em ii}) of~\cite{BS97a}.
There, we also show that $\eta\preceq\zeta$ if and only if
\begin{enumerate}
\item $a\in$up$(\zeta)\Longrightarrow\eta(a)\leq\zeta(a)$.
\item $b\in$down$(\zeta)\Longrightarrow\eta(b)\geq\zeta(b)$.
\item $a,b\in$up$(\zeta)$ (or
$a,b\in$down$(\zeta)$) with $a<b$ and
$\zeta(a)<\zeta(b)\Longrightarrow\eta(a)<\eta(b)$.
\end{enumerate}

\begin{lem}\label{lem:perm_facts}
Let $\zeta\in{\mathcal S}_\infty$.
The labeled poset $[1,\zeta]_\preceq$ is a chain if and only if 
$\zeta$ is either an increasing or a decreasing cycle.
Moreover, if $\zeta$ is an increasing (decreasing) cycle of length $m+1$,
then the chain 
$[1,\zeta]_\preceq$ is increasing (decreasing) and $\zeta$ is
shape-equivalent to $r[m,1]$  ($v(1^m,m)$).
\end{lem}

\noindent{\bf Proof. }
Let $\zeta\in{\mathcal S}_\infty$ and construct 
$u\leq_k\zeta u$ as above.
Set $m:=\ell(\zeta u)-\ell(u)$, and consider any chain in $[u,w]_k$:
$$
u=u_0\stackrel{b_1}{\longrightarrow} u_1
\stackrel{b_2}{\longrightarrow} u_2
\ \cdots\  u_{m-1}\stackrel{b_m}{\relbar\joinrel\longrightarrow} u_m=w.
$$

Suppose that the poset $[1,\zeta]_\preceq\simeq[u,\zeta u]_k$ is a chain.
By Order 2, $\zeta$ is irreducible.
We show that $\zeta$ is either an increasing or a decreasing cycle 
by induction on $m$.
Suppose $\eta=u_{m-1}u^{-1}$ is an increasing cycle.
Then $\eta=(b_{m-1},b_{m-2},\ldots,b_1,a_1)$ where
$u_1=(a_1,b_1)u$ and $u_i=(b_{i-1},b_i)u_{i-1}$ for $i>1$.
Let $\zeta=(a_m,b_m)\eta$.

Since $u_{m-1}^{-1}(b_{m-1})\leq k$ and $u_{m-1}^{-1}(b_m)> k$,
we must have $b_{m-1}\neq b_m$.
If $b_m>b_{m-1}$ so that $[1,\zeta]_\preceq$ is increasing,
then, as $\zeta$ is irreducible, we must have $a_m=b_{m-1}$ and so
$\zeta$ is the increasing cycle
$$
(b_m,b_{m-1},\ldots,b_1,a_1).
$$
Indeed, if either $a_m>b_{m-1}$ or $a_m<b_{m-2}$, then 
$[1,\zeta]_\preceq$ is not a chain, and
$b_{m-1}>a_m\geq b_{m-2}$ contradicts
$u_{m-2}\lessdot_k u_{m-1}\lessdot_k u_m$.
Suppose now that $b_m<b_{m-1}$, then the irreducibility of $\zeta$
implies that $m=2$ and $b_m=a_1$, so that 
$[1,\zeta]_\preceq$ is decreasing and $\zeta$ is a decreasing cycle.

Similar arguments suffice when $\eta=u_{m-1}u^{-1}$ is a decreasing cycle,
and the other statements are straightforward.
\QED

\noindent{\bf Proof that Identities 1 and 2 imply the Pieri formula. }

Let $\zeta\in{\mathcal S}_\infty$ 
and suppose $c^\zeta_{(m,0,\ldots,0)}\neq 0$.
Then $m=|\zeta|$, by homogeneity.
Replacing $\zeta$ by a shape-equivalent permutation if necessary, we 
may assume that $\zeta\in{\mathcal S}_n$ and
$\zeta(i)\neq i$ for each $1\leq i\leq n$.

Define $u$ and $w:=\zeta u$ as in the first paragraph of this section, so
that $u,w\in{\mathcal S}_n$ and 
$c^\zeta_{(m,0,\ldots,0)}=c^w_{u\,r[m,k]}$.
Since $c^w_{u\,r[m,k]}\neq 0$, we must have
$m=n-k=\#$down$(\zeta)$: 
Consider any chain
\begin{equation}\label{eq:beta-chain}
u=u_0\stackrel{b_1}{\longrightarrow} u_1
\stackrel{b_2}{\longrightarrow} u_2
\ \cdots \ u_{m-1}\stackrel{b_m}{\relbar\joinrel\longrightarrow} u_m=w
\end{equation}
in $[u,w]_k$.
Then down$(\zeta)\subset\{b_1,\ldots,b_m\}$ so that
$m\geq n-k$.
However, $c^w_{u\,r[m,k]}\neq 0$ and $w\in{\mathcal S}_n$
implies that  
$r[m,k]\in{\mathcal S}_n$, and hence 
$k{+}m\leq n$.
It follows that down$(\zeta)=\{b_1,\ldots,b_m\}$.
Thus if we have $u_i=u_{i-1}(c_i,\,d_i)$ with $c_i\leq k<d_i$, then 
by the construction of $u$, 
$\{d_1,\ldots,d_m\} = \{k{+}1,\ldots,k{+}m=n\}$.

Consider the case when $\zeta$ is irreducible.
Then we must have $c_1=c_2=\cdots=c_m$.
This implies that $k=\#$up$(\zeta) = 1$, and $m=n-1$.
By (1) of Definition~\ref{def:1} we must then have 
$b_1<b_2<\cdots<b_m$, and hence
$\zeta=(n,\,n{-}1,\,\ldots,\,2,\,1)$, an increasing cycle.
But this is $r[n{-}1,1]$, so $u=1$, the identity permutation.
Since $c^w_{1\,v} = \delta_{w,\,v}$, the Kronecker delta, 
$c^\zeta_{\lambda}=\delta_{\lambda,\,(m,0,\ldots,0)}$ and so
$S_\zeta = h_{n-1}$.

If more generally we have $\eta\in{\mathcal S}_\infty$ with
$\#$down$(\eta)=|\eta|=m$ and $\eta$ irreducible, 
then considering a shape-equivalent $\zeta\in{\mathcal S}_n$ with 
$n$ minimal, we see that $\eta$ is an increasing cycle and 
$S_{\eta}=h_m$.

We return to the case of 
$\zeta\in{\mathcal S}_n$ with $c^\zeta_{(m,0,\ldots,0)}\neq 0$.
Let $\zeta=\zeta_1\cdots\zeta_t$ be the disjoint 
factorization of $\zeta$ into irreducibles.
Then each $\zeta_i$ is an increasing cycle.
Suppose that $m_i=|\zeta_i|$.
By Identity 2, we have that 
\begin{eqnarray*}
S_\zeta&=& S_{\zeta_1}\cdots S_{\zeta_t}\\
&=& h_{m_1}\cdots h_{m_t}.
\end{eqnarray*}
This is equivalent to~\cite[Theorem 5]{Sottile96}.
{}From this, we deduce that $c^\zeta_\lambda = c^\mu_{\nu\,\lambda}$,
where $\mu/\nu$ is a horizontal strip with $m_i$ boxes in the 
$i$th row.
By the classical Pieri formula for Schur polynomials, this implies that 
$c^\zeta_{(m,0,\ldots,0)}=1$.
\QED

\section{Skew Schur functions from labeled posets}

In~\cite[Theorem 4.3]{BS_monoid}, we showed how the Pieri formula implies
Identity 1. 
Here we complete the proof of Theorem~\ref{thm:main_equiv}, showing how the
Pieri formula implies Identity 2.
The first step is a reinterpretation of a construction 
in~\cite[\S4]{BS_monoid} from which we associate a symmetric function
to any symmetric labeled poset.
For intervals in Young's lattice, we obtain skew Schur functions, and for
intervals in either a $k$-Bruhat order or the $\preceq$-order, skew Schubert
functions.
In Section 4, we show that for intervals in the weak order we obtain
Stanley symmetric functions.

Let $P$ be a labeled poset with total rank $m$.
A (maximal) chain in $P$ gives a sequence of edge labels, called the 
{\em word} of that chain.
A {\em composition}  $\alpha:= (\alpha_1,\ldots,\alpha_k)$ of 
$m=\alpha_1+\cdots+\alpha_k$ ($\alpha_i\geq 0$), 
determines, and is determined by a (multi)subset 
$I(\alpha):=\{\alpha_1,\alpha_1+\alpha_2,\ldots,\alpha_1+\cdots+\alpha_k\}$
of $\{1,\ldots,m\}$.
For a composition  $\alpha$  of $m=$ rank$P$, let 
$H_{\alpha}(P)$ be the set of (maximal) chains in $P$ 
whose word $w$ has descent set $\{j\mid w_j>w_{j+1}\}$ 
contained in the set $I(\alpha)$.
We adopt the convention that the last position of a word is a descent.
If some $\alpha_i<0$, then we set $H_{\alpha}(P)=\emptyset$.
We say that $P$ is {\em (label-) symmetric} if the cardinality of
$H_{\alpha}(P)$ depends only upon the parts of 
$\alpha$ and not their order.

Let $\Lambda$ be the ${\mathbb Z}$-algebra of symmetric functions.
Recall that $\Lambda = {\mathbb Z}[h_1,h_2,\ldots]$, where $h_i$ is the
complete homogeneous symmetric function of degree $i$, the sum of all
monomials of degree $i$.
For a composition $\alpha$, set 
$$
h_{\alpha}\ :=\ h_{\alpha_1}h_{\alpha_2}\cdots h_{\alpha_k}.
$$ 

\begin{definition}
{\sl
Suppose $P$ is a symmetric labeled poset.
Define the ${\mathbb Z}$-linear map 
$\chi_P : \Lambda \rightarrow {\mathbb Z}$ by
$$
\chi_P\ :\ h_{\alpha} \longmapsto \# (H_{\alpha}(P)).
$$
For any partition $\lambda$, define the skew coefficient $c^P_\lambda$ to be
$\chi_P(S_\lambda)$, where $S_\lambda$ is the Schur symmetric function.
}
\end{definition}

We point out some properties of these coefficients $c^P_\lambda$.
For a partition $\lambda$ of $m$ ($\lambda\vdash m$) with $\lambda_{k+1}=0$
and a permutation $\pi\in {\mathcal S}_k$, let $\lambda_\pi$ be the
following composition of $m$:
$$
\pi(1)-1+\lambda_{k+1-\pi(1)},\,\pi(2)-2+\lambda_{k+1-\pi(2)},\,
\ldots,\,\pi(k)-k+\lambda_{k+1-\pi(k)}.
$$

\begin{lem}\label{lem:coeff}
Let $P,Q$ be symmetric labeled posets.
\begin{enumerate}
\item For any partition $\lambda$,
$$
   c^P_\lambda\ :=\ \sum_{\pi\in {\mathcal S}_k}
   \varepsilon(\pi) \#( H_{\lambda_\pi}(P))
$$
where $\lambda_{k+1}=0$ and 
$\varepsilon :{\mathcal S}_k \to \{\pm 1\}$ is the sign character.
\item If $P\simeq Q$ as labeled posets
(Definition~\ref{def:lposet}) then for any partition $\lambda$, 
$c^P_\lambda = c^Q_\lambda$.
\end{enumerate}
\end{lem}

The first statement follows from the Jacobi-Trudi formula, and the second
by noting that the bijection $P\leftrightarrow Q$ induces bijections
$H_{\alpha}(P)\leftrightarrow H_{\alpha}(Q)$.

\begin{rem}
By the Pieri formula for Schubert polynomials, the number
$\#(H_{\alpha}([u,w]_k))$ is the coefficient of 
${\mathfrak S}_w$ in the product
${\mathfrak S}_u\cdot h_{\alpha}(x_1,\ldots,x_k)$.
It follows that intervals in a $k$-Bruhat order or in the $\preceq$-order
are symmetric.
For similar reasons, we see that intervals in Young's lattice are symmetric,
as $\#(H_{\alpha}([\mu,\lambda]_\subset))$  is the skew Kostka coefficient
$K_{\alpha,\,\lambda/\mu}$, which is 
the coefficient of $S_\lambda$ 
in $S_\mu\cdot h_\alpha$, equivalently, the number of semistandard Young
tableaux of shape $\lambda/\mu$ and content $\alpha$.
One may construct an explicit bijection with the second set as follows:
A chain in $H_{\alpha}([\mu,\lambda]_\subset)$ is naturally decomposed
into subchains with increasing labels of lengths 
$\alpha_1,\alpha_2,\ldots,\alpha_k$.
Placing the integer $i$ in the boxes corresponding to covers in 
the $i$th such subchain furnishes the bijection.
\end{rem}

\begin{prop}[Theorem~4.3 of~\cite{BS_monoid}]
Let $u\leq_k w$ and $\lambda\vdash \ell(w)-\ell(u)= m$.
Then $c^w_{u\,v(\lambda,k)} = c^{[u,w]_k}_\lambda$.
\end{prop}

\noindent{\bf Proof. }
By definition, $c^w_{u\,v(\lambda,k)}$ is the 
coefficient of ${\mathfrak S}_w$ in
the expansion of the product
${\mathfrak S}_u\cdot S_\lambda(x_1,\ldots,x_k)$ into Schubert polynomials.
By the Jacobi-Trudi formula,
\begin{eqnarray*}
{\mathfrak S}_u\cdot S_\lambda(x_1,\ldots,x_k) &=&
{\mathfrak S}_u\cdot \sum_{\pi\in {\mathcal S}_k} 
\varepsilon(\pi) h_{\lambda_\pi}(x_1,\ldots,x_k) \\
 &=& \sum_w \sum_{\pi\in {\mathcal S}_k} 
\varepsilon(\pi) \# (H_{\lambda_\pi}([u,w]_k))\: {\mathfrak S}_w\\
 &=& \sum_w  c^{[u,w]_k}_\lambda\; {\mathfrak S}_w.
\qquad \QED
\end{eqnarray*}

\begin{prop}[Corollary 4.9 of\/~\cite{BS_monoid}]\label{prop:identity}
If $u\leq_kw$ and $y\leq_l z$ with $wu^{-1}$ shape equivalent to 
$zy^{-1}$, then for all $\lambda$,
$c^w_{u\,v(\lambda,k)} = c^z_{y\,v(\lambda,l)}$.
\end{prop}

\noindent{\bf Proof. }
By Order 1, $[u,w]_k\simeq [y,z]_l$ is an isomorphism of labeled posets.
\QED

\begin{defn}\label{def:skew}
Let $P$ be a ranked labeled poset with total rank $m$.
Define the symmetric function $S_P$ by
$$
S_P\ :=\ \sum_{\lambda\vdash m} c^P_\lambda S_\lambda,
$$
where $S_\lambda$ is a Schur {\em function}.
\end{defn}

\noindent{\bf Proof of Theorem~\ref{thm:skew} (1), (2), and (3). }
(1) is a consequence of Lemma~\ref{lem:coeff} (2).
For (3), let $\mu\subset \nu$ in Young's lattice, suppose $\nu_{k+1}=0$, and 
consider the interval $[\mu,\nu]_\subset$ in Young's lattice.
Then $[\mu,\nu]\simeq[v(\mu,k),v(\nu,k)]_k$, and so 
$c^{[\mu,\nu]}_\lambda = c^{v(\nu,k)}_{v(\mu,k)\,v(\lambda,k)}=
c^\nu_{\mu\,\lambda}$.
Hence $S_{[\mu,\nu]_\subset}=S_{\nu/\mu}$.
Similarly, we see that for $u\leq_k w$ or $\zeta\in{\mathcal S}_\infty$,
we have $S_{[u,w]_k}=S_{wu^{-1}}$ and $S_{[1,\zeta]_\preceq}=S_\zeta$, 
the skew Schubert functions of \S 1.
\QED

\begin{rem}\label{rem:cyclic}
According to Proposition~\ref{prop:identity}, the skew Schubert function
$S_\zeta$ depends only on the shape equivalence class of $\zeta$.
In~\cite{BS97a} there is another identity:\medskip

Theorem H of\/~\cite{BS97a}.
{\em 
Suppose $\eta,\zeta\in{\mathcal S}_n$ with $\zeta=\eta^{(12\ldots n)}$.
Then   $S_\eta=S_\zeta$.
}\medskip

The example of $\eta=(1243)$ and $\zeta=(1243)$ in ${\mathcal S}_4$
(see Figure~\ref{fig:interval})
shows that in general 
$[1,\eta]_\preceq \not\simeq [1,\eta^{(12\ldots n)}]_\preceq$.
However, these two intervals do have the same number of maximal
chains~\cite[Corollary 1.4]{BS97a}.
In fact, for $\eta\in{\mathcal S}_n$ and $\alpha$ a composition,
$\# (H_{\alpha}([1,\eta]_\preceq))=
\# (H_{\alpha}([1,\eta^{(12\ldots n)}]_\preceq))$.

Thus if $\sim$ is the equivalence relation generated by shape equivalence
and this `cyclic shift' 
($\eta\sim\eta^{(12\ldots n)}$, if $\eta\in{\mathcal S}_n$), 
then $S_\zeta$ depends only upon the $\sim$-equivalence class of $\zeta$.
(This is analogous to, but stronger than the fact that the skew Schur
function  $S_\kappa$ depends on $\kappa$ only up to a translation in the
plane.)

There is a combinatorial object $\Gamma_\zeta$ which determines
the $\sim$-equivalence class of $\zeta$.
First place the set $\{a\mid a\neq \zeta(a)\}$ at the vertices of a regular
$\#\{a\mid a\neq \zeta(a)\}$-gon in clockwise order.
Next, for each $a$ with $a\neq \zeta(a)$,  draw a directed chord from $a$ 
to $\zeta(a)$.
$\Gamma_\zeta$ is the resulting configuration of directed chords, 
up to rotation and dilation and without any vertices labeled
({\em cf.}~\cite[\S3.3]{BS97a}).
The irreducible factors of $\zeta$ correspond to connected 
components of $\Gamma_\zeta$ (considered as a subset of the plane).
The figure $\Gamma_{(1243)}=\Gamma_{(1423)}$ is also displayed in 
Figure~\ref{fig:interval}.

\begin{figure}[htb]
    $$\epsfxsize=3.5in \epsfbox{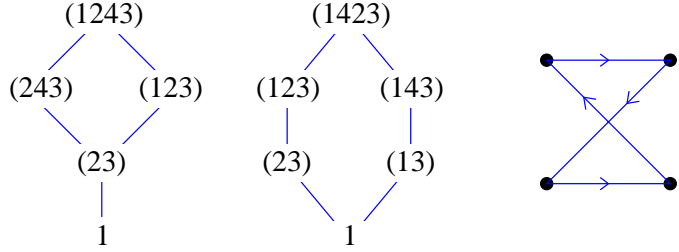}$$
    \caption{Intervals under cyclic shift and 
    $\Gamma_\zeta$\label{fig:interval}}
\end{figure}
\end{rem}

We conclude this section with the following Theorem:

\begin{thm}\label{thm:product}
Let $P$ and $Q$ be symmetric labeled posets with disjoint sets of 
edge labels.
Then
$$
S_{P\times Q}\ =\ S_P\cdot S_Q.
$$
\end{thm}

This will complete the proof of Theorem~\ref{thm:main_equiv}, namely that 
the Pieri formula and Order 2 imply Identity 2:
If $\zeta\cdot\eta$ is a disjoint product, then 
$[1,\zeta]_\preceq$ and $[1,\eta]_\preceq$ have disjoint sets of edge
labels.
Together with Theorem~\ref{thm:skew}(4), this gives another proof of
Theorem~3.4 in~\cite{Stanley84}, that $F_{w\times u}=F_w\cdot F_u$.
\medskip

To prove Theorem~\ref{thm:product}, we first study chains in 
$H_{\alpha}(P\times Q)$.
Suppose that $P$ has rank $n$ and $Q$ has rank $m$.
Note that a chain in $P\times Q$ determines and is determined by the
following data:
\begin{equation}\label{chain:data}
 \begin{array}{l}
  \bullet\ \mbox{A chain in each of $P$ and $Q$},\\
  \bullet\ \mbox{A subset $B$ of $\{1,\ldots,n+m\}$ with $\#B=n$}.
 \end{array}
\end{equation}
Recall that covers $(p,q)\lessdot(p',q')$ in $P\times Q$ 
have one of two forms:
either $p=p'$ and $q'$ covers $q$ in $Q$ or else
$q=q'$ and $p'$ covers $p$ in $P$.
Thus a chain in $P\times Q$ gives a chain in each of $P$ and $Q$,
with the covers from $P$ interspersed among the covers from $Q$.
If we set $B$ to be the positions of the covers from $P$, we 
obtain the description~(\ref{chain:data}).
Define 
$$
\mbox{sort}\ :\ \mbox{\em chains}(P\times Q)\ 
\longrightarrow\  \mbox{\em chains}(P)\times\mbox{\em chains}(Q)
$$
to be the map which forgets the positions $B$ of the covers from 
$P$.

\begin{lem}\label{lem:bijection}
Let $P$ and $Q$ be labeled posets with disjoint sets of edge labels
and $\alpha$ be any composition.
Then
$$
\mbox{\rm sort}\ :\ H_{\alpha}(P\times Q)\ \longrightarrow\ 
\coprod_{\beta+\gamma=\alpha}
H_{\beta}(P)\times H_{\gamma}(Q)
$$
is a bijection.
\end{lem}

For integers $a<b$, let $[a,b]:=\{n\in{\mathbb Z}\mid a\leq n\leq b\}$.
For a chain $\xi$, let $\xi|_{[a,b]}$ be the
portion of $\xi$ starting at the $a$th step and continuing to the 
$b$th step.

\noindent{\bf Proof. }
Let $\xi\in H_{\alpha}(P\times Q)$
and set $I=I(\alpha)$ so that $I_i=\alpha_1+\cdots+\alpha_i$.
Then sort$(\xi)\in H_{\beta}(P)\times H_{\gamma}(Q)$,
where, for each $i$,  $\beta_i$ counts
the number of covers of $\xi|_{[I_{i-1},I_i]}$
from $P$ and $\gamma_i=\alpha_i-\beta_i$.

To see this is a bijection, we construct its inverse.
For chains $\xi^P\in H_{\beta}(P)$ and
$\xi^Q\in H_{\gamma}(Q)$ with
$\beta+\gamma=\alpha$,  define the set $B$ by the
conditions 
\begin{enumerate}
\item $\beta_i=\# B\cap [I(\alpha)_{i-1},I(\alpha)_i]$.
\item If $b_1\leq\cdots\leq b_{\beta_i}$ and
$c_1\leq\cdots\leq c_{\gamma_i}$
are the covers in 
$\xi^P|_{[I(\beta)_{i-1},I(\beta)_i]}$ and
$\xi^Q|_{[I(\gamma)_{i-1},I(\gamma)_i]}$ respectively,
then, up to a shift of $I(\alpha)_{i-1}$, the set
$B\cap [I(\alpha)_{i-1},I(\alpha)_i]$ records the positions of the 
the $b$'s in the linear ordering of
$\{b_1,\ldots,b_{\beta_1},c_1,\ldots,c_{\gamma_i}\}$.
\end{enumerate}
This clearly gives the inverse to the map sort.
\QED

Recall that the comultiplication 
$\Delta:\Lambda \rightarrow \Lambda\otimes\Lambda$ is defined by
$$
\Delta(h_a)\ =\ \sum_{b+c=a}h_b\otimes h_c.
$$
Thus, for a composition $\alpha$, 
$$
\Delta(h_\alpha)\ =\ \sum_{\beta+\gamma=\alpha}h_\beta\otimes h_\gamma.
$$

{}From Lemma~\ref{lem:bijection}, we immediately deduce:
\begin{cor}
Let $P,Q$ be symmetric labeled posets with disjoint sets of edge labels.
Then
$$
\begin{picture}(100,67)
\put( 0,25){$\chi_{P\times Q}$}
\put(9,43){$\Lambda$}
\put(43,0){${\mathbb Z}$}
\put(43,52){$\Delta$}
\put(77,43){$\Lambda\otimes\Lambda$}
\put(72,25){$\chi_P\otimes \chi_Q$}
\put(18,40){\vector(1,-1){27}}
\put(80,40){\vector(-1,-1){27}}
\put(19,47){\vector(1,0){55}}
\end{picture}
$$
commutes.
\end{cor}

\begin{cor}
Let $P,Q$ be symmetric labeled posets with disjoint sets of edge labels.
Then, for any partition $\lambda$, 
$$
c^{P\times Q}_\lambda\ =\ 
\sum_{\mu,\nu} c^\lambda_{\mu\,\nu}\: c^P_\mu\; c^Q_\nu.
$$
\end{cor}

\noindent
{\bf Proof. }
Recall~\cite[I.5.9]{Macdonald95} that
$\Delta(S_\lambda)\ =\ \sum_{\mu,\nu} c^\lambda_{\mu\,\nu}\:S_\mu\;S_\nu$.
Hence
$$
\chi_{P\times Q}(S_\lambda)\ =\ 
\sum_{\mu,\nu} c^\lambda_{\mu\,\nu}\: \chi_P(S_\mu)\:\chi_P(S_\nu).\qquad
\QED
$$

We complete the proof of Theorem~\ref{thm:product}:
Let $P,Q$ be symmetric labeled posets with disjoint sets of edge labels.
Then
\begin{eqnarray*}
S_P\cdot S_Q&=&
\sum_{\mu,\nu}c^P_\mu\;S_\mu\: c^Q_\nu\; S_\nu\\
&=& \sum_{\lambda,\mu,\nu}c^\lambda_{\mu\,\nu}\;c^P_\mu\;c^Q_\nu\:S_\lambda\\
&=&\sum_\lambda c^{P\times Q}_\lambda S_{\lambda}\ \ =\ \ 
S_{P\times Q}.\qquad
\QED
\end{eqnarray*}

\section{Stanley symmetric functions from labeled posets}

We establish Theorem~\ref{thm:skew}(4)
by adapting the proof of the Littlewood-Richardson rule 
in~\cite{Remmel_Shimozono} to obtain a bijective interpretation of the 
constants $c^{[1,w]_{\mbox{\scriptsize weak}}}_\lambda$,
which shows $S_{[1,w]_{\mbox{\scriptsize weak}}}=F_w$ by the formulas 
in~\cite{LS82b,EG}.
The main tool is a {\em jeu de taquin} for reduced decompositions. 

We use Cartesian conventions for Young diagrams and 
skew diagrams.
Thus the first row is at the bottom.
A filling of a diagram $D$ with positive integers which increase across 
rows and up columns is a {\em tableau} with {\em shape} $D$.
The {\em word} of a tableau is the sequence of its entries, read across 
each row starting with the topmost row.

A {\em reduced decomposition} $\rho$ for a permutation 
$w\in{\mathcal S}_\infty$ is the 
word of a maximal chain in 
$[1,w]_{\mbox{\scriptsize\rm weak}}$.
Let $R(w)$ be the set of all reduced decompositions for $w$
and for a composition $\alpha$ of $\ell(w)$, write $H_\alpha(w)$ 
for $H_\alpha([1,w]_{\mbox{\scriptsize\rm weak}})$.
Given any composition $\alpha$ and any 
reduced decomposition $\rho\in H_{\alpha}(w)$, there is a unique smallest
diagram $\lambda/\mu$ with row lengths 
$\lambda_i-\mu_i=\alpha_{k+1-i}$ for which 
$\rho$ is the word of a tableau $T(\alpha,\rho)$ 
of shape $\lambda/\mu$.
By this we mean that $\mu_j-\mu_{j+1}$ is minimal for all $j$.
If $\mu_1=0$, then $T(\alpha,\rho)$ has {\em partition shape}
$\lambda$ ($=\alpha$), otherwise  $T(\alpha,\rho)$ has {\em skew shape}.
Given a reduced decomposition $\rho\in R(w)$, define $T(\rho)$ to be the 
tableau $T(\alpha,\rho)$, where 
$I(\alpha)$ is the descent set of $\rho$.

Stanley~\cite{Stanley84} defined a symmetric function $F_w$ 
for every $w\in{\mathcal S}_\infty$.
(That $F_w$ is symmetric includes a proof that the intervals 
$[1,w]_{\mbox{\scriptsize weak}}$ are symmetric.)
Thus there exists integers $a^w_\lambda$ such that
$$
F_w\ =\ \sum_{\lambda\vdash l} a^w_\lambda S_\lambda.
$$
A combinatorial interpretation for $a^w_\lambda$
was given (independently) in~\cite{LS82b} and~\cite{EG}:
$$
a^w_\lambda\ =\ \#\{ \rho\in R(w)\mid T(\rho)\mbox{ has 
partition shape }\lambda\}.
$$
(See~\cite[\S VII]{Macdonald91} for an account with proofs.)
Theorem~\ref{thm:skew}(4) is a consequence of the following result:

\begin{thm}\label{thm:weak_coefficients}
For any $w\in{\mathcal S}_\infty$ and partition $\lambda\vdash\ell(w)$,
$$
a^w_\lambda \ =\ c^{[1,w]_{\mbox{\scriptsize\rm weak}}}_\lambda.
$$
\end{thm}

Our proof is based on the proof of the 
Littlewood-Richardson rule given by Remmel and 
Shimozono~\cite{Remmel_Shimozono}.
We define an involution $\theta$ on the set
$$
\coprod_{\pi\in{\mathcal S}_k} \{\pi\}\times 
H_{\lambda_\pi}(w)
$$
(here $\lambda\vdash \ell(w)$ and $\lambda_{k+1}=0$)
such that 
\begin{enumerate}
\item
$\theta(\pi,\rho) = (\pi,\rho)$ if and only if 
$T(\rho)$ has shape $\lambda$, from which it follows that $\pi=1$.
\item If $T(\rho)$ does not have shape $\lambda$, then
$\theta(\pi,\rho) = (\pi',\rho')$
where $T(\rho')$ does not have shape $\lambda$ and
$\rho'\in H_{\lambda_{\pi'}}(w)$ with $|\ell(\pi)-\ell(\pi')|=1$.
\end{enumerate}

Theorem~\ref{thm:weak_coefficients} is 
a corollary of the existence of such an involution $\theta$: 
By property 2, only the fixed points 
of $\theta$ contribute to the sum in 
Lemma~\ref{lem:coeff}(1).

The involution $\theta$ will be defined using a {\em jeu de taquin} for 
tableaux whose words are reduced decompositions.
Because we only play this {\it jeu de taquin} on diagrams with 
two rows, we do not describe it in full.

\begin{defn}
Let $T$ be a tableau of shape $(y+p,q)/(y,0)$ whose word is a reduced 
decomposition for a permutation $w$.
If $y\neq 0$, we may perform an inward slide.
This modification of an ordinary {\it jeu de taquin}
slide ensures we obtain a tableau whose 
word is a reduced decomposition of $w$.

Begin with an empty box at position $(y,1)$ and move it through the tableau
$T$ according to the following local rules:
\begin{enumerate}
\item If the box is in the first row, it switches with whichever 
of its neighbors to the right or above is smaller.

If both neighbors are equal, say they are $a$, then their other 
neighbor is necessarily $a+1$, as we have a reduced decomposition.
Locally we will have the following configuration,
where \raisebox{-2pt}{%
\begin{picture}(10,10)  
\put( 0,10){\line(0,-1){10}}\put(10,10){\line(-1,0){10}}
\put( 0, 0){\line(1, 0){10}}\put(10, 0){\line( 0,1){10}}
\put( 0, 0){\line(1, 1){10}}\put(10, 0){\line(-1,1){10}}
\end{picture}}
denotes the empty box and $a+b+1<c$:
$$
\begin{picture}(200,30)   \thicklines
\put( 0, 0){\line(1, 1){15}}\put(15, 0){\line(-1,1){15}}
\put( 0, 0){\line(1,0){200}}  \put( 0,15){\line(1,0){200}}
\put( 0,30){\line(1,0){184}}
\put( 0, 0){\line(0,1){30}}   \put(15, 0){\line(0,1){30}}
\put(42, 0){\line(0,1){30}}   \put(72, 0){\line(0,1){30}}
\put(114,15){\line(0,1){15}}  \put(142, 0){\line(0,1){30}}
\put(184, 0){\line(0,1){30}}  \put(200, 0){\line(0,1){15}}
\put(4,19){$a$}
\put(26,4){$a$}               \put(19,19){$a{+}1$}
\put(46,4){$a{+}1$}           \put(46,19){$a{+}2$}         
\put(86,4){$\cdots\cdots$}    \put(86,19){$\cdots$}
\put(118,19){$a{+}b$}
\put(153,4){$a{+}b$}          \put(146,19){$a{+}b{+}1$}
\put(188,4){$c$}
\end{picture}
$$
The empty box moves through this configuration, transforming it into:
$$
\begin{picture}(200,30)
\thicklines
\put( 0, 0){\line(1,0){198}}  \put( 0,15){\line(1,0){198}}
\put( 0,30){\line(1,0){184}}
\put( 0, 0){\line(0,1){30}}   \put(27, 0){\line(0,1){30}}
\put(56, 0){\line(0,1){30}}   
\put(100,0){\line(0,1){30}}   \put(142, 0){\line(0,1){30}}
\put(184, 0){\line(0,1){30}}  \put(198, 0){\line(0,1){15}}
\put(142,15.5){\line(3, 1){42}}\put(142,29.5){\line(3,-1){42}}
\put( 9, 4){$a$}              \put( 4,19){$a{+}1$}
\put(31,4){$a{+}1$}           \put(31,19){$a{+}2$}         
\put(63,4){$\cdots\cdots$}    \put(63,19){$\cdots\cdots$}
\put(111, 4){$a{+}b$}         \put(104,19){$a{+}b{+}1$}
\put(146,4){$a{+}b{+}1$}
\put(188,4){$c$}
\end{picture}
$$
This guarantees that we still have a reduced decomposition for $w$.
\item If the box is in the second row, then it switches with its neighbour 
to the right.

\end{enumerate}

If $y+p>q$, then we may analogously perform an outward slide, 
beginning with an empty box at $(q+1,2)$ and sliding to the 
left or down according to local rules that are the reverse of  those for the
inward slide. 
\end{defn}

We note some consequences of this definition.
\begin{itemize}
\item
The box will change rows at the first pair of entries $b\leq c$ 
it encounters with $b$ at $(i,2)$ and $c$ 
immediately to its lower right at $(i+1,1)$.
If there is no such pair, it will change rows at the end of 
the first row in an inward slide if $p+y=q$, and at the beginning 
of the second row in an outward slide if $y=0$.

\item 
At least one of these will occur if $y$ is minimal given 
the word of the tableau and $p,q$.
Suppose this is the case.
Then the tableau $T'$ obtained from a slide will 
have another such pair $b'\leq c'$ with $b'$ at $(\imath',2)$ and $c'$ at 
$(\imath'+1,1)$.
Hence, if we perform a second slide, the box will again change rows.

\item The inward and outward slides are inverses.

\end{itemize}

Let $\overline{H}_\alpha(w)$ be the subset of $H_\alpha(w)$ consisting 
of chains $\rho$ such that $T(\alpha,\rho)$ has skew shape.
The proof of the following lemma is straightforward.

\begin{lem}\label{lem:bijections}
Let $w\in {\mathcal S}_\infty$
and suppose $p<q$ with $p+q=\ell(w)$.
Then $H_{(q,p)}(w)=\overline{H}_{(q,p)}(w)$
and
\begin{enumerate}
\item 
For every  
$\rho\in H_{(q,p)}(w)$, we may perform $q-p$ 
inward slides to $T((q,p),\rho)$.
If $\rho'$ is the word of the resulting tableau, then the map 
$\rho\mapsto \rho'$ defines a bijection
$$
H_{(q,p)}(w)\ 
\longleftrightarrow \ 
H_{(p,q)}(w).
$$
The inverse map is given by the application of $q-p$ outward slides.
\item
If we now let $\rho'$ be the word of the tableau obtained
after $q-p-1$ inward slides to $T((q,p),\rho)$ for 
$\rho\in H_{(q,p)}(w)$, then the map 
$\rho\mapsto \rho'$ defines a bijection
$$
\overline{H}_{(q,p)}(w)\ 
\longleftrightarrow \ 
\overline{H}_{(p+1,q-1)}(w).
$$
The inverse map is defined by the application of $q-1-p$ outward slides.
\end{enumerate}
\end{lem}

The first part gives a proof that intervals in the weak
order are symmetric:
Let $\alpha=(\alpha_1,\ldots,\alpha_k)$ and 
$\alpha'=(\alpha_1,\ldots,\alpha_{r+1},\alpha_r,\ldots,\alpha_k)$ 
be compositions of $\ell(w)$.
Then applying the bijection in Lemma~{lem:bijections}(1) to the segment
$\rho_r$ of $\rho\in H_{\alpha}(w)$
between $I(\alpha)_{r-1}$ and  $I(\alpha)_{r+1}$
defines a bijection
$$
H_{\alpha}(w)\ \longleftrightarrow \ H_{\alpha'}(w).
$$

\begin{rem}
This bijection is different from the one used 
in~\cite{Stanley84} to prove symmetry of these intervals.
Indeed, consider the example given there,
which we write as a tableau:
$$
\epsfxsize=3.5in \epsfbox{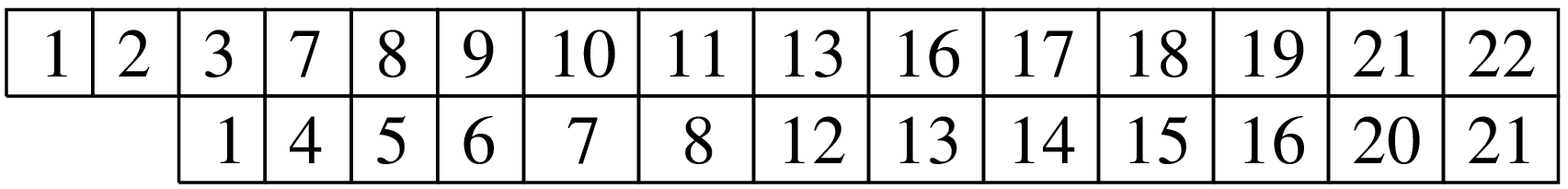}
$$
In~\cite{Stanley84}, Stanley maps this to 
$$
\epsfxsize=4.03in \epsfbox{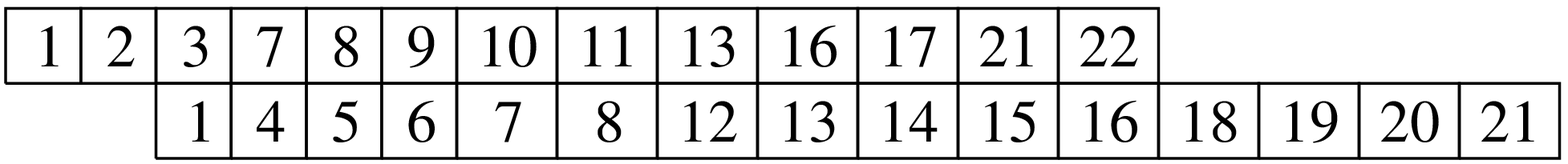}
$$
But the bijection we define gives us this:
$$
\epsfxsize=3.63in \epsfbox{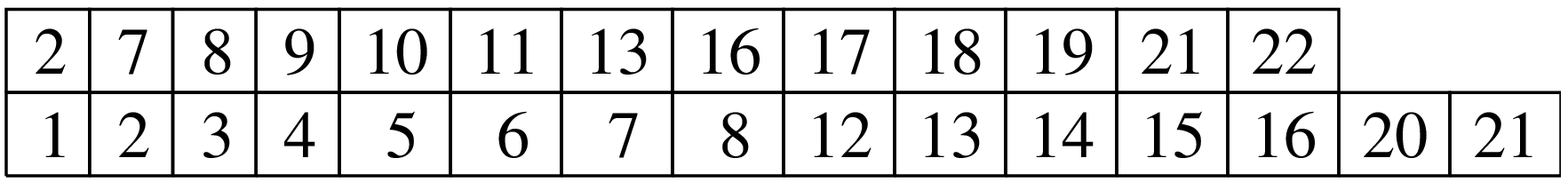}
$$
\end{rem}
Now we may define $\theta$.
By the definition of $\lambda_\pi$,
if $\rho\in H_{\lambda_\pi}(w)$, then $T(\rho)$ has shape $\lambda$ if 
and only if $T(\lambda_\pi,\rho)$ has partition shape, which implies that 
$\pi=1$.

\begin{defn}\label{def:theta}
Suppose $w\in{\mathcal S}_\infty$ and $\lambda\vdash \ell(w)$ is a partition
with $\lambda_{k+1}=0$.
Let $\pi\in{\mathcal S}_k$.
For $\rho\in H_{\lambda_\pi}(w)$,
define $\theta(\pi,\rho)$ as follows:
\begin{enumerate}
\item If $T(\rho)$ has shape $\lambda$, set
$\theta(\pi,\rho)=(\pi,\rho)$.
In this case, $\pi=1$, so  $\lambda_\pi=\lambda$ and
$T(\rho)=T(\lambda_\pi,\rho)$.
\item If $T(\rho)$ does not have shape $\lambda$,
then $T(\lambda_\pi,\rho)$ has skew shape and we select 
$r=r(T(\lambda_\pi,\rho))$ with $1\leq r<k$ as follows:\smallskip

Left justify the rows of $T(\lambda_\pi,\rho)$.
Since $T(\lambda_\pi,\rho)$ has skew shape, there is an entry $a$ of this 
left-justified figure in postiton $(i,r+1)$ either with no entry in 
position $(i,r)$ just below it, or else with an entry $b\geq a$ just below 
it.
Among all such $(i,r)$ choose the one with $i$ minimal, and for this $i$,
$r$ maximal. 
\smallskip

Let $\rho_r$ be the word given by the rows $r+1$ and $r$ of
$T(\lambda_\pi,\rho)$ and $(q,p)$ the lengths of these two rows.
Then $T((q,p),\rho_r)$ has skew shape, and we may apply the map of 
Lemma~\ref{lem:bijections}(2) to obtain the word $\rho_r'$.
Define $\theta(\pi,\rho)=(\pi',\rho')$, where
$\rho'$ is the word obtained from $\rho$ by replacing $\rho_r$ 
with $\rho_r'$ and $\pi'\pi^{-1}=(r,\,r{+}1)$.
Note that $T(\lambda_{\pi'},\rho')$ also has skew shape and 
$T(\rho')$ does not have shape $\lambda$.
\end{enumerate}

\end{defn}

\begin{ex}
Let $w=4621357$ and $\lambda=(4,3,3,1)$.
Then 
$\rho=5.345.236.1236\in H_{\lambda}(w)$	
but 
$$
\begin{picture}(135,60)
\thicklines
\put( 0,22.5){$T(\lambda,\rho)\ =$}
\put( 60,30){\line(0,1){30}}   \put(60,60){\line(1,0){15}}
\put( 75, 0){\line(0,1){60}}   \put(60,45){\line(1,0){45}}
\put( 90, 0){\line(0,1){45}}   \put(60,30){\line(1,0){60}}
\put(105, 0){\line(0,1){45}}   \put(75,15){\line(1,0){60}}
\put(120, 0){\line(0,1){30}}   \put(75, 0){\line(1,0){60}}
\put(135, 0){\line(0,1){15}}
\put(64,49){5}
\put(64,34){3}\put(79,34){4}\put(94,34){5}
              \put(79,19){2}\put(94,19){3}\put(109,19){6}
              \put(79, 4){1}\put(94, 4){2}\put(109, 4){3}\put(124, 4){5}
\end{picture}
$$
has skew shape.
Left-justifying the rows of $T(\lambda,\rho)$, we obtain:
$$
\begin{picture}(60,60)
\thicklines
\put( 0, 0){\line(0,1){60}}   \put(0,60){\line(1,0){15}}
\put(15, 0){\line(0,1){60}}   \put(0,45){\line(1,0){45}}
\put(30, 0){\line(0,1){45}}   \put(0,30){\line(1,0){45}}
\put(45, 0){\line(0,1){45}}   \put(0,15){\line(1,0){60}}
\put(60, 0){\line(0,1){15}}   \put(0, 0){\line(1,0){60}}
\put(4,49){5}
\put(4,34){3}\put(19,34){4}\put(34,34){5}
\put(4,19){2}\put(19,19){3}\put(34,19){6}
\put(4, 4){1}\put(19, 4){2}\put(34, 4){3}\put(49, 4){5}
\end{picture}
$$
This is not a tableau, as the third column reads $365$, which is not 
increasing.
Since this is the first such column and the last decrease
is at position $2$, we have $r=2$.
Since these two rows each have length 3, we perform 
one outward slide (by our choice of $r$, we can perform 
such a slide!) to obtain the tableau $T((4,2),\rho'_r)$ as follows:
$$
\begin{picture}(60,30)
\thicklines
\put( 0,15){\line(0,1){15}}   \put(15, 0){\line(1,0){45}}
\put(15, 0){\line(0,1){30}}   \put( 0,15){\line(1,0){60}}
\put(30, 0){\line(0,1){30}}   \put( 0,30){\line(1,0){60}}
\put(45, 0){\line(0,1){30}}  
\put(60, 0){\line(0,1){30}}  
\put(45,15){\line(1,1){15}}   \put(45,30){\line(1,-1){15}}  
\put(4,19){3}\put(19,19){4}\put(34,19){5}
             \put(19, 4){2}\put(34, 4){3}\put(49, 4){6}
\end{picture}
\qquad \raisebox{14pt}{$\relbar\joinrel\longrightarrow$} \qquad
\begin{picture}(60,30)
\thicklines
\put( 0,15){\line(0,1){15}}   \put(15, 0){\line(1,0){45}}
\put(15, 0){\line(0,1){30}}   \put( 0,15){\line(1,0){60}}
\put(30, 0){\line(0,1){30}}   \put( 0,30){\line(1,0){60}}
\put(45, 0){\line(0,1){30}}  
\put(60, 0){\line(0,1){30}}  
\put(45, 0){\line(1,1){15}}   \put(45,15){\line(1,-1){15}}  
\put(4,19){3}\put(19,19){4}\put(34,19){5}\put(49,19){6}
             \put(19, 4){2}\put(34, 4){3}
\end{picture}
\qquad \raisebox{14pt}{$\relbar\joinrel\longrightarrow$} \qquad
\begin{picture}(45,30)
\thicklines
\put( 0,15){\line(0,1){15}}   \put(15, 0){\line(1,0){45}}
\put(15, 0){\line(0,1){30}}   \put( 0,15){\line(1,0){60}}
\put(30, 0){\line(0,1){30}}   \put( 0,30){\line(1,0){60}}
\put(45, 0){\line(0,1){30}}  
\put(60, 0){\line(0,1){30}}  
\put(15, 0){\line(1,1){15}}   \put(15,15){\line(1,-1){15}}  
\put(4,19){3}\put(19,19){4}\put(34,19){5}\put(49,19){6}
                           \put(34, 4){2}\put(49, 4){3}
\end{picture}
$$

Thus $\rho'=5.3456.23.1235\in
H_{\lambda_{(2,\,3)}}(w)$.
If we left justify $T(\lambda_{(2,\,3)},\rho')$, then 
we obtain:
$$
\begin{picture}(60,60)
\thicklines
\put( 0, 0){\line(0,1){60}}   \put(0,60){\line(1,0){15}}
\put(15, 0){\line(0,1){60}}   \put(0,45){\line(1,0){60}}
\put(30, 0){\line(0,1){45}}   \put(0,30){\line(1,0){60}}
\put(45, 0){\line(0,1){15}}   \put(0,15){\line(1,0){60}}
\put(60, 0){\line(0,1){15}}   \put(0, 0){\line(1,0){60}}
\put(45,30){\line(0,1){15}}  
\put(60,30){\line(0,1){15}}   
\put(4,49){5}
\put(4,34){3}\put(19,34){4}\put(34,34){5}\put(49,34){6}
\put(4,19){2}\put(19,19){3}
\put(4, 4){1}\put(19, 4){2}\put(34, 4){3}\put(49, 4){5}
\end{picture}
$$
The 5 in the third row has no lower neighbour, hence
$2=r(\lambda,\rho)=r(\lambda_{(2,\,3)},\rho')$.
\end{ex}

We complete the proof of Theorem~\ref{thm:weak_coefficients}
by showing that $\theta$ is an involution.
This is a consequence of Lemma~\ref{lem:bijections}(2) and the following fact:

\begin{lem}
In  (2) of Definition~\ref{def:theta}, 
if $\rho\in H_{\lambda_\pi}(w)$
and $T(\lambda_{\pi},\rho)$ has skew shape, then 
$r(T(\lambda_{\pi},\rho))=r(T(\lambda_{\pi'},\rho'))$.
\end{lem}

\noindent{\bf Proof. }
Suppose we are in the situation of (2) in Definition~\ref{def:theta}.
The lemma follows once we show that that 
$T((q,p),\rho_r)$ and $T((p+1,q-1),\rho'_r)$ agree in the first 
$i$ entries of their second rows, the first $i-1$ entries of their
first rows, and the 
$i$th entry $c$ in the first row of $T((p+1,q-1),\rho'_r)$ satisfies
$a\leq c$, or else there is no $i$th entry.

In fact, we show this holds for each intermediate tableau
obtained from $T((q,p),\rho)$ by some of the slides used to form
$T((p+1,q-1),\rho')$.

We argue in the case that $p<q$, that is, for inward slides.
Suppose that $T$ is an intermediate tableau satisfying the claim, and 
that the tableau $T'$ obtained from $T$ by a single inward slide is 
also an intermediate tableau.
It follows that $T'$ has skew shape, so that if $(y+s,t)/(y,0)$ is the shape
of $T$, then $y>1$.

Suppose that during the slide the box changes rows
at the $j$th column.
We claim that $j\geq i+y-1(>i)$.
If this occurs, then the first $i$  entries in the second row and 
first $i-1$ entries in the first row of $T$ are unchanged in $T'$.
Also, the $i$th entry in the first row of $T'$  is either the $i$th
entry in the first row of $T$ (if $j\geq i+y$) or it is the $j$th
entry in the second row of $T$, which is greater than the $i$th
entry, $a$.
Thus showing  $j\geq i+y-1$ completes the proof.

To see that $j\geq i+y-1$ note that if $j$ is the last column, then 
$j=t=s+y$.
Since $s\geq i-1$, we see that $j\geq y+i-1$.
If $j$ is not the last column, then the entries $b$ at $(j,2)$ 
$c$ at $(j+1,1)$ of $T$ satisfy $b\leq c$.
Suppose that $j<i+y-1$.
Then $c$ is the ($j-y+1$)th entry in the first row of 
$T$.   
Since $j-y+1<i$, our choice of $i$ ensures that $c$ is less than the 
entry at $(j-y+1,2)$ of $T$.
Since $j-y+1<j$, this in turn is less than $b$, a contradiction.

Similar arguments suffice for the case when $p\geq q$.
\QED

\begin{rem}
While it may seem this proof has only a formal relation to 
the proof of Remmel and Shimozono~\cite{Remmel_Shimozono},
it is in fact nearly an exact translation---the only difference
being in our choice of $r$.
(Their choice of $r$ is not easily expressed in this setting.)
We elaborate.

The exact same proof, but with the ordinary 
{\em jeu de taquin}, shows that 
$c^{[\mu,\lambda]_\subset}_\nu$ counts the chains in 
$[\mu,\lambda]_\subset$ whose word is the word of a tableau of shape $\nu$.
This is just the Littlewood-Richardson coefficient $c^\lambda_{\mu\,\nu}$.
One way to see this is to consider the bijection between 
$H_\nu([\mu,\lambda]_\subset)$ and the set of semistandard Young tableaux
of shape $\lambda/\mu$ and content $(\nu_k,\ldots,\nu_1)$.
The chains whose word is the word of a tableau of shape $\nu$
correspond to {\em reverse} LR tableaux of shape $\lambda/\mu$, which
are defined as follows:

Let $f_{a,b}(T)$ be the number of $a$'s in the first $b$ positions of 
the word of $T$.
A reverse LR tableau $T$ with largest entry $k$ is a tableau satisfying:
$$
f_{1,b}(T)\ \leq\ f_{2,b}(T)\ \leq\ \cdots\ \leq\ f_{k,b}(T)
$$
for all $b$.
It is an exercise to verify that there are exactly
$c^\lambda_{\mu\nu}$ reverse LR tableaux of shape $\lambda/\mu$ and 
content $\nu_k,\ldots,\nu_2,\nu_1$.

The choice we make of $i$ and $r$ is easily expressed in these terms:
$i$ is the minimum value of $f_{a,b}(T)$ among all violations
$f_{a,b}(T)>f_{a+1,b}(T)$, and if $a$ is the minimal first index among
all violations with $f_{a,b}(T)=i$, then $r=k-a$.
The choice in~\cite{Remmel_Shimozono} for reverse LR tableaux
would be $r=k-a$, where $f_{a,b}(T)$ is the violation with minimal $b$.

The key step we used was the {\em jeu de taquin}
whereas Remmel and Shimozono used an operation built from the 
$r$-pairing of Lascoux and Sch\"utzenberger~\cite{LS81}.
In fact, this too is a direct translation.

The reason for this is, roughly, that the passage from the word of a 
chain $\rho\in H_\alpha([\mu,\lambda]_\subset)$ to a semistandard Young
tableau of shape $\lambda/\mu$ and content $(\alpha_k,\ldots,\alpha_1)$ 
(which interchanges shape with content) also interchanges Knuth 
equivalence and dual Knuth equivalence~\cite{Haiman_dual_equivalence}.
The operators constructed from the 
$r$-pairing preserve the dual equivalence class of a 2-letter word but 
alter its content.
In fact, this property characterizes such an operation.

As shown in~\cite{Haiman_dual_equivalence}, there is at most one tableau
in a given Knuth equivalence class and a given dual equivalence class.
Also, for semistandard Young tableaux with at most 2 letters, there is 
at most one tableau with given partition shape and content.
It follows that any operation on tableaux acting on the 
subtableau of entries $r,r+1$ which preserves the dual 
equivalence class of the subtableau, but reverses its content
is uniquely defined by these properties.

Thus the symmetrization operators in~\cite{LS81},
which generate an ${\mathcal S}_\infty$-action on tableaux extending 
the natural action on their contents, is unique.
Expressed in this form, we see that this action coincides with 
one introduced earlier by Knuth~\cite{Knuth}.
This action was the effect of permuting rows of a matrix
on the $P$-symbol obtained from that matrix by Knuth's generalization of 
the Robinson-Schensted correspondence.
The origin of these symmetrization operators in the work of Knuth has been
overlooked by most authors, perhaps because Bender-Knuth~\cite{Bender_Knuth}
later use a different operation to prove symmetry.
\end{rem}

For each poset $P$ in the classes of labeled posets we consider 
here, the symmetric function $S_P$ is Schur-positive.
When $P$ is an interval in some $k$-Bruhat order, this follows from
geometry, for intervals in Young's lattice, this is a consequence of the 
Littlewood-Richardson rule, and for intervals in the weak order, it is 
due to Lascoux-Sch\"utzenberger~\cite{LS82a} and Edelman-Greene~\cite{EG}.
Is there a representation-theoretic explanation?
In particular, we ask:\medskip 

\noindent{\bf Question:}
{\em 
If $P$ is an interval in a $k$-Bruhat order, can one construct a 
representation $V_P$ of ${\mathcal S}_{\mbox{\scriptsize\rm rank}P}$
so that $S_P$ is its Frobenius character?
More generally, for a labeled poset $P$, can one define a (virtual) 
representation $V_P$ so that $S_P$ is its Frobenius character?
If so, is $V_{P\times Q}\simeq V_P\otimes V_Q$?
}\medskip

When P is an interval in Young's lattice this is a skew Specht module.
For an interval $[1,w]_{\mbox{\scriptsize\rm weak}}$ in the weak order,
Kr\'askiewicz~\cite{Kraskiewicz} constructs a 
${\mathcal S}_{\ell(w)}$-representation of dimension $\# R(w)$.
For general linear group represenations, such a construction is known.
For intervals in the weak oder, this is due to 
Kr\'askiewicz and Pragacz~\cite{KP}.

\section{The monomials in a Schubert polynomial}

We give a new proof based upon geometry that a Schubert polynomial is a 
sum of monomials with non-negative coefficients.
This analysis leads to a combinatorial construction of Schubert 
polynomials in terms of chains in the Bruhat order.
It also shows these coefficients are certain intersection numbers,
essentially the same interpretation found by Kirillov and Maeno~\cite{KM}.

The first step is Theorem~\ref{thm:univariate}, which generalizes both 
Proposition 1.7 of~\cite{LS82b} and 
Theorem C ({\em ii}) of~\cite{BS97a}.
Recall that 
$u\stackrel{r[m,k]}{\relbar\joinrel\relbar\joinrel\longrightarrow} w$
when one of the following equivalent conditions holds:
\begin{itemize}
\item $c^w_{u,\,r[m,k]}=1$.
\item $u\leq_k w$ and $wu^{-1}$ is a disjoint product of increasing
cycles.
\item There is an chain in $[u,w]_k$:
$$
u\ \stackrel{b_1}{\longrightarrow}\ u_1
\stackrel{b_2}{\longrightarrow}\ \cdots\ 
\stackrel{b_m}{\longrightarrow}\ u_m=w
$$
with $b_1<b_2<\cdots<b_m$.
\end{itemize}

For $p\in{\mathbb N}$, define the map 
$\Phi_p:{\mathbb Z}[x_1,x_2,\ldots]\longrightarrow
{\mathbb Z}[y]\otimes {\mathbb Z}[x_1,x_2,\ldots]$
by
$$
\Phi_p(x_i)\ =\ \left\{\begin{array}{ll}
x_i    &\ \mbox{if}\ i<p\\
y      &\ \mbox{if}\ i=p\\
x_{i-1}&\ \mbox{if}\ i>p   \end{array}\right..
$$
For $w\in{\mathcal S}_\infty$ and $p,q\in{\mathbb N}$,
define $\varphi_{p,q}(w)\in{\mathcal S}_\infty$ by
$$
\varphi_{p,q} (w)(j) \quad =\quad \left\{\begin{array}{lcl}
w(j)     && j < p \mbox{ and } w(j)  <  q\\
w(j)+1   && j < p \mbox{ and } w(j)\geq q\\
q        && j = p\\
w(j-1)   && j > p \mbox{ and } w(j)  <  q\\
w(j-1)+1 && j > p \mbox{ and } w(j)\geq q
\end{array}\right..
$$
Representing permutations as matrices, $\varphi_{p,q}$ adds a
new $p$th row and $q$th column consisting mostly of zeroes, but with a 1
in the $(p,q)$th position.
For example,
$$
\varphi_{3,3}(23154)\ =\ 243165 \qquad\mbox{and}\qquad
\varphi_{2,5}(2341)\ =\ 25342. 
$$

\begin{thm}\label{thm:univariate}
For $u\in{\mathcal S}_n$,
$$
\Phi_p {\mathfrak S}_u\ =\ 
\sum_{\stackrel{\mbox{\scriptsize $j, w$ \rm with}}%
{u\stackrel{r[n{+}1{-}p{-}j,\,p]}%
{\relbar\joinrel\relbar\joinrel\relbar\joinrel\relbar\joinrel%
\relbar\joinrel\relbar\joinrel\relbar\joinrel\relbar\joinrel\longrightarrow}
\varphi_{p,n+1}(w)}}
y^j\:{\mathfrak S}_w(x).
$$
Moreover, if $n$ is not among $\{u(1),\ldots,u(p-1)\}$, then 
the sum may be taken over those $j,w$ with
$u\stackrel{r[n{-}p{-}j,\,p]}{\relbar\joinrel\relbar\joinrel%
\relbar\joinrel\relbar\joinrel\relbar\joinrel\longrightarrow}
\varphi_{p,n}(w)$.
\end{thm}

Iterating this gives another proof that the
monomials in a Schubert polynomial have non-negative coefficients.

\begin{ex}
Consider $\Phi_2{\mathfrak S}_{13542}$.
We display all increasing chains in the 2-Bruhat order on
${\mathcal S}_5$ above
$13542$ whose endpoint $w$ satisfies $w(2)=5$:
$$\epsfxsize=2.in \epsfbox{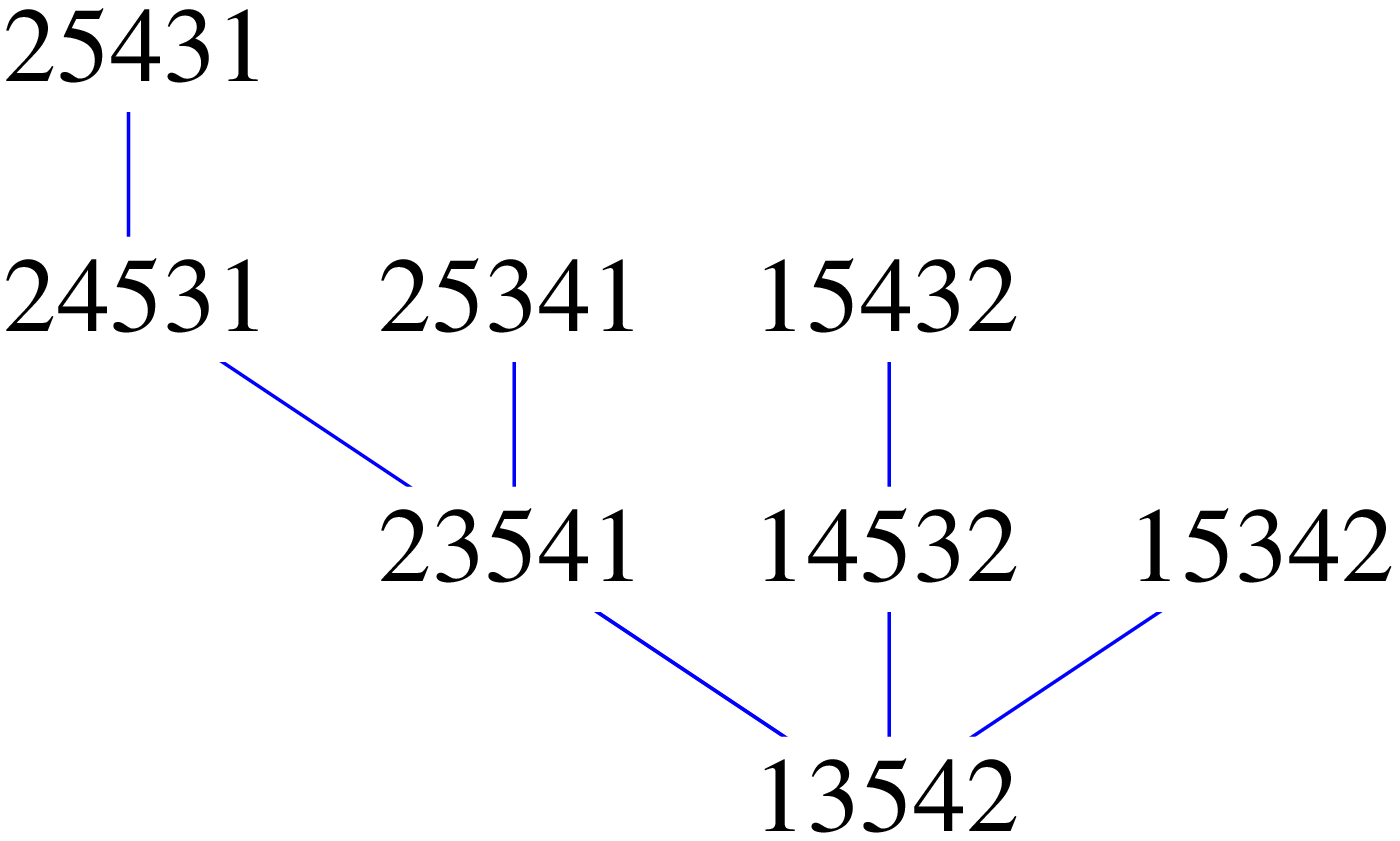}$$
We see therefore that 
\begin{eqnarray*}
13542\ \stackrel{r[3,2]}{\relbar\joinrel\relbar\joinrel\longrightarrow}
25431&= & \varphi_{2,5}(2431),\\
13542\ \stackrel{r[2,2]}{\relbar\joinrel\relbar\joinrel\longrightarrow}
25341&= & \varphi_{2,5}(2341),\\
13542\ \stackrel{r[2,2]}{\relbar\joinrel\relbar\joinrel\longrightarrow}
15432&= & \varphi_{2,5}(1432),\\
13542\ \stackrel{r[1,2]}{\relbar\joinrel\relbar\joinrel\longrightarrow}
15342&= & \varphi_{2,5}(1342).
\end{eqnarray*}
Then Theorem~\ref{thm:univariate}
asserts that
$$
\Phi_2{\mathfrak S}_{13542}\ =\ 
{\mathfrak S}_{2431}(x) + y{\mathfrak S}_{2341}(x) +
y{\mathfrak S}_{1432}(x) + y^2{\mathfrak S}_{1342}(x),
$$
which may also be verified by direct calculation.
\end{ex}

\noindent{\bf Proof of Theorem~\ref{thm:univariate}. }
We make two definitions.
For $p\leq n$, define another map 
$\psi_{p,[n]}:{\mathcal S}_n\times {\mathcal S}_m
\hookrightarrow {\mathcal S}_{n+m}$ by
\begin{equation}\label{eq:psi-map}
  \psi_{p,[n]}(w,z)(i)\ =\ \left\{\begin{array}{ll}
  w(i)    &\ i<p\\
  n+z(1)  &\ i=p\\
  w(i-1)  &\ p<i\leq n+1\\
  n+z(i-n)&\ n+1<i\leq n+m   \end{array}\right..
\end{equation}
Then $\psi_{p,[n]}(1,1)=r[n{+}1{-}p,p]$.

Let $P\subset\{1,2,\ldots,n+m\}$ and suppose that
\begin{eqnarray*}
P&=&p_1<p_2<\cdots<p_n,\\
\{1,\ldots,n+m\}- P&=& q_1<q_2<\cdots<q_m.
\end{eqnarray*}
Define the map 
$\Psi_P:{\mathbb Z}[x_1,x_2,\ldots,x_{n+m}]\longrightarrow
{\mathbb Z}[x_1,\ldots,x_n]\otimes{\mathbb Z}[y_1,\ldots,y_m]$
by
$$
\Psi_P(x_i)\ =\ \left\{\begin{array}{rl}
x_j&\ \mbox{if}\ i=p_j\\
y_j&\ \mbox{if}\ i=q_j
\end{array}\right..
$$

Suppose now that $P=\{1,2,\ldots,p-1,p+1,\ldots,n+1\}$.
Then  for $u\in{\mathcal S}_{n+m}$,
Theorem 4.5.4 of~\cite{BS97a} asserts that

\begin{equation}\label{eq:coh}
\Psi_P{\mathfrak S}_u\ \equiv\ 
\sum_{w\in{\mathcal S}_n,\ z\in{\mathcal S}_m}
c^{\psi_{p,[n]}(w,z)}_{u\ r[n{+}1{-}p,p]}
{\mathfrak S}_w(x)\otimes{\mathfrak S}_z(y),
\end{equation}
modulo the ideal 
$\langle {\mathfrak S}_w(x)\otimes 1, 1\otimes{\mathfrak S}_z(y)
\mid w\not\in{\mathcal S}_n,z\not\in{\mathcal S}_m\rangle$
which is equal to the ideal 
$\langle x^{\alpha}\otimes 1, 1\otimes y^{\alpha}
\mid \alpha_i\geq n-i\mbox{ for some }i\rangle$.
(The calculation is in the cohomology of the product of flag manifolds
{\em Flags}$({\mathbb C}\,^n)\times${\em Flags}$({\mathbb C}\,^m)$.)

Suppose now that $u\in{\mathcal S}_n$ and $m\geq n$.
Then~(\ref{eq:coh}) is an identity of polynomials, and not just of
cohomology classes.
We also see that $\Psi_P{\mathfrak S}_u=\Phi_p{\mathfrak S}_u$, since
${\mathfrak S}_u\in{\mathbb Z}[x_1,\ldots,x_n]$.
By the Pieri formula,
$$
c^{\psi_{p,[n]}(w,z)}_{u\ r[n{+}1{-}p,p]}\ =\ 
\left\{\rule{0pt}{20pt}\right.
\begin{array}{ll} 
1 & \ \mbox{if }u\stackrel{r[n{+}1{-}p,\,p]}{\relbar\joinrel\relbar\joinrel%
\relbar\joinrel\relbar\joinrel\relbar\joinrel\longrightarrow}
\psi_{p,[n]}(w,z),\\
0&\ \mbox{otherwise}.\rule{0pt}{12pt}
\end{array}
$$
Since $u\leq_p \psi_{p,[n]}(w,z)$ and $u(n+i)=n+i$, 
Definition~\ref{def:1} (2) (for $u\leq_p \psi_{p,[n]}(w,z)$)
implies that
$$
\psi_{p,[n]}(w,z)(n+1)<\psi_{p,[n]}(w,z)(n+2)<\cdots.
$$
Thus by the definition~(\ref{eq:psi-map}) of $\psi_{p,[n]}$, we have 
$z(2)<z(3)<\cdots$, and so $z$ is the Grassmannian
permutation $r[z(1){-}1,1]$.
Hence ${\mathfrak S}_z(y)=y^{z(1)-1}$.

If we set $j=z(1)-1$, then 
$\psi_{P,[n]}(w,z)=\varphi_{p,n+1+j}(w)$.
Thus, for $u\in{\mathcal S}_n$, we have 
$$
\Phi_p{\mathfrak S}_u\ =\ 
\sum_{\stackrel{\mbox{\scriptsize $j,w$ such that}}%
{u\stackrel{r[n{+}1{-}p,\,p]}%
{\relbar\joinrel\relbar\joinrel\relbar\joinrel\relbar\joinrel%
\relbar\joinrel\relbar\joinrel\longrightarrow}
\varphi_{p,n+1+j}(w)}}
y^j\,{\mathfrak S}_w(x).
$$
Suppose that 
$u\stackrel{r[n{+}1{-}p,\,p]}%
{\relbar\joinrel\relbar\joinrel\relbar\joinrel\relbar\joinrel%
\relbar\joinrel\longrightarrow}
\varphi_{p,n+1+j}(w)$.
Consider the unique increasing chain in the interval
$[u,\ \varphi_{p,n+1+j}(w)]_p$:
$$
u=u_0\stackrel{b_1}{\longrightarrow}\ \cdots\ 
\stackrel{b_{n-p-j}}{\relbar\joinrel\relbar\joinrel%
\relbar\joinrel\longrightarrow}
u_{n-p-j}\stackrel{b_{n+1-p-j}}{\relbar\joinrel\relbar\joinrel%
\relbar\joinrel\relbar\joinrel\longrightarrow}
\ \cdots\ \stackrel{b_{n+1-p}}{\relbar\joinrel\relbar%
\joinrel\relbar\joinrel\longrightarrow}
\varphi_{p,n+1+j}(w).
$$
Because $u\in{\mathcal S}_n$, we must have $b_{n+1-p-j}=n+1$
and so $u_{n+1-p-j}=\varphi_{p,n+1}(w)$.
Moreover, if $n$ is not among $\{u(1),\ldots,u(p)\}$, then 
we have  $b_{n-p-j}=n$
and so $u_{n-p-j}=\varphi_{p,n}(w)$.
If $u(p)=n$, then we also have $u_{n-p-j}=\varphi_{p,n}(w)$.
This completes the proof.
\QED

Define $\delta$ to be the sequence $(n-1,n-2,\ldots,1,0)$.

\begin{cor}\label{cor:chain_monomial}
For $w\in{\mathcal S}_n$ and $\alpha<\delta$, the coefficient of
$x^{\delta-\alpha}$ in ${\mathfrak S}_w$ is the number of chains
$$
  w\lessdot w_1\lessdot w_2\lessdot\cdots\lessdot
  w_{\alpha_1+\cdots+\alpha_{n-1}} = \omega_0
$$
in the Bruhat order where, for each $1\leq k\leq n-1$,
\begin{equation}\label{eq:ch-cond}
  w_{\alpha_1+\cdots+\alpha_{k-1}}\:\lessdot_k\: 
w_{1+\alpha_1+\cdots+\alpha_{k-1}}\:
  \lessdot_k\:\cdots\:\lessdot_k\: w_{\alpha_1+\cdots+\alpha_k} 
\end{equation}
is an increasing chain in the $k$-Bruhat order.
\end{cor}

\begin{ex}
Here are all such chains in ${\mathcal S}_4$ from $1432$ to 
$4321$, with the index $\alpha$ displayed above each chain:
$$\epsfxsize=2.in \epsfbox{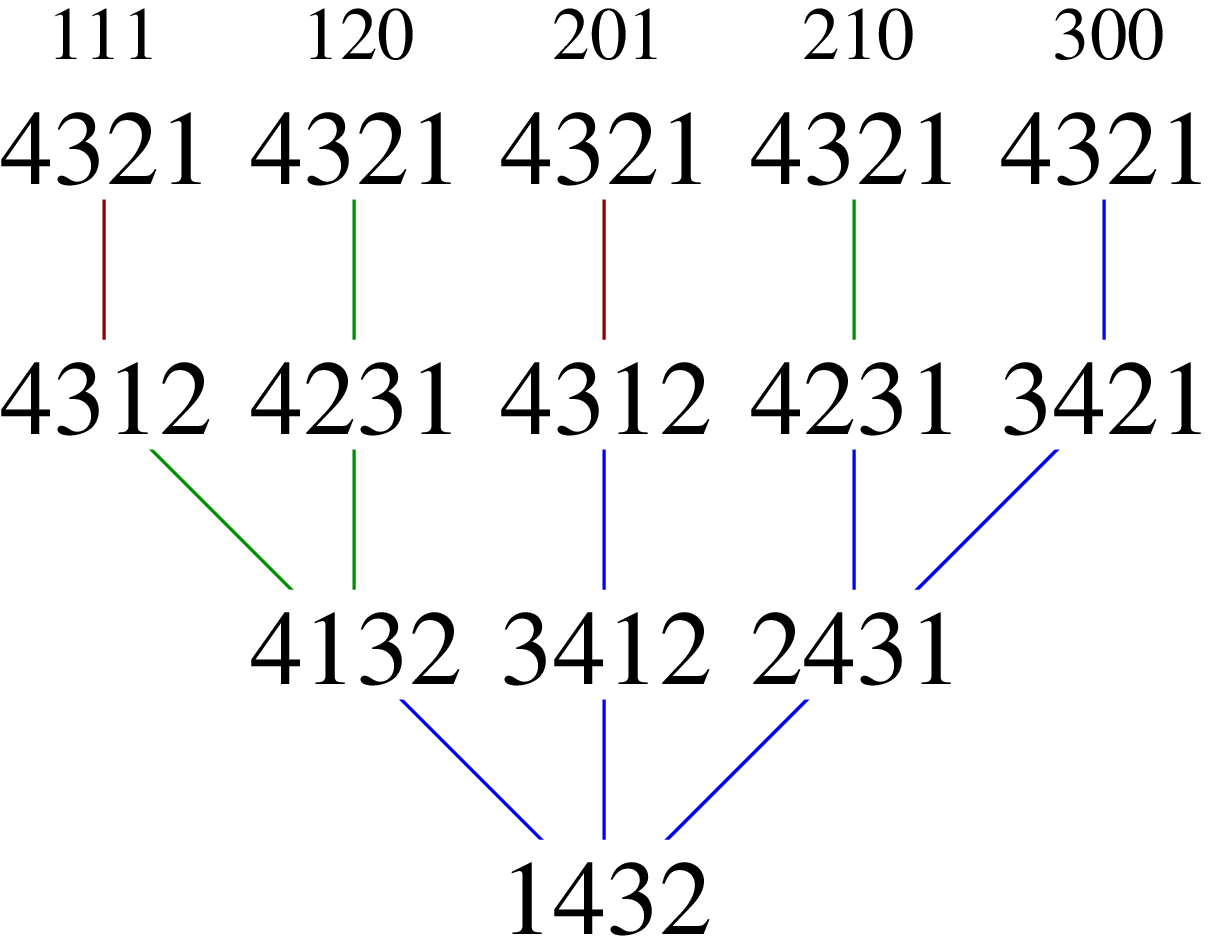}$$

{}From this, we see that 
\begin{eqnarray*}
{\mathfrak S}_{1432}&=&
x^{321-111}+x^{321-120}+x^{321-201}+x^{321-210}+x^{321-300}\\
&=&x_1^2x_2 + x_1^2x_3 + x_1x_2^2 + x_1x_2x_3 + x_2^2x_3.
\end{eqnarray*}
\end{ex}

\noindent{\bf Proof. }
Repeatedly applying $\Phi_1$ and iterating Theorem~\ref{thm:univariate},
we see that the coefficient of $x^{\delta-\alpha}$ in
${\mathfrak S}_w(x)$ is the number of chains
$$		
  w\lessdot w_1\lessdot w_2\lessdot\cdots\lessdot
  w_{\alpha_1+\cdots+\alpha_{n-1}} = \omega_0
$$		
which satisfy the conditions of the corollary, together with the 
(apparent) additional requirement that, for each $k<n$,
\begin{equation}
\label{eq:value}
w_{\alpha_1+\cdots+\alpha_k}(j)\ =\ n+1-j\ \mbox{ for all }\ j \leq k.
\end{equation}
The corollary will follow, once we show this is no additional restriction.

First note that if 
$u\stackrel{r[a,k]}{\relbar\joinrel\relbar\joinrel\longrightarrow}u'$
with $u'(j)=n+1-j$ for $1\leq j\leq k$, but
$u(i)<n+1-i$ for some $1\leq i\leq k$, then $i=k$.
To see this, note that since $u\leq_k u'$, the form of 
$u'$ and Definition~\ref{def:1} (2)
implies that $u(1)>u(2)>\cdots>u(k)$.
Set $\zeta=u'u^{-1}$.
Since $u\stackrel{r[a,k]}{\relbar\joinrel\relbar\joinrel\longrightarrow}u'$,
$\zeta$ is a disjoint product of increasing cycles,
hence their supports are are non-crossing.
Suppose $i<k$. 
Then $\{u(i),n+1-i=u'(i)\}$ and 
$\{u(i+1),n-i=u'(i+1)\}$
are in the support of distinct cycles.
However, $u(i+1)<u(i)\leq n-i<n+1-i$ contradicts that these supports are 
non-crossing, so we must have $i=k$.

Let 
$$
  w\lessdot w_1\lessdot w_2\lessdot\cdots\lessdot
  w_{\alpha_1+\cdots+\alpha_{n-1}} = \omega_0
$$
be a chain which satisfies the conditions of the corollary.
We prove that~(\ref{eq:value}) holds for all $k<n$ by downward induction.
Since $\omega_0=w_{\alpha_1+\cdots+\alpha_{n-1}}$, 
we see that~(\ref{eq:value}) holds 
for $k=n-1$.
Suppose that~(\ref{eq:value}) holds for some $k$.
Set $u=w_{\alpha_1+\cdots+\alpha_{k-1}}$ and 
$u'=w_{\alpha_1+\cdots+\alpha_k}$.
Then 
$u\stackrel{r[\alpha_k,k]}{\relbar\joinrel\relbar\joinrel\longrightarrow}u'$
with $u'(j)=n+1-j$ for $1\leq j\leq k$.
By the previous paragraph, we must have $u(i)=n+1-i$ for 
all $i<k$, hence~(\ref{eq:value}) holds for $k-1$.
\QED

We could also have written the coefficient of $x^{\delta-\alpha}$
in ${\mathfrak S}_w(x)$ as the number of chains
$$
w\ \stackrel{r[\alpha_1,1]}{\relbar\joinrel%
\relbar\joinrel\relbar\joinrel\longrightarrow}\ 
w_1\ \stackrel{r[\alpha_2,2]}{\relbar\joinrel%
\relbar\joinrel\relbar\joinrel\longrightarrow}\ 
w_2\ \stackrel{r[\alpha_3,3]}{\relbar\joinrel%
\relbar\joinrel\relbar\joinrel\longrightarrow}\ 
\cdots\  \stackrel{r[\alpha_{n-1},n{-}1]}{\relbar\joinrel\relbar\joinrel%
\relbar\joinrel\relbar\joinrel\relbar\joinrel\longrightarrow}\ 
\omega_0
$$
in ${\mathcal S}_n$.
{}From this and the Pieri formula for Schubert polynomials,
we obtain another description of these coefficients.
First, for $\alpha=(\alpha_1,\alpha_2,\ldots,\alpha_{n-1})$ with 
$\alpha_i\geq 0$,
let $h(\alpha)$ denote the product of complete homogeneous
symmetric polynomials
$$
h_{\alpha_1}(x_1) h_{\alpha_2}(x_1,x_2)\cdots
h_{\alpha_{n-1}}(x_1,x_2,\ldots,x_{n-1}).
$$

\begin{cor}
For $w\in{\mathcal S}_n$,
$$
{\mathfrak S}_w \ =\ \sum_{\alpha}
d^w_{\alpha}x^{\delta-\alpha}
$$
where $d^w_{\alpha}$ is the coefficient of 
${\mathfrak S}_{\omega_0}$ in the product
${\mathfrak S}_w\cdot h(\alpha)$.
\end{cor}

This is essentially the same formula as found by Kirillov and
Maeno~\cite{KM} who showed that the coefficient of 
$x^{\delta-\alpha}$ in ${\mathfrak S}_w$ to be the coefficient
of ${\mathfrak S}_{\omega_0}$ in the product
${\mathfrak S}_{\omega_0 w\omega_0}\cdot e(\alpha)$, where
$$
e(\alpha)=e_{\alpha_{n-1}}(x_1)e_{\alpha_{n-2}}(x_1,x_2)\cdots
e_{\alpha_1}(x_1,\ldots,x_{n-1}).
$$
To see these are equivalent, note that the algebra involution
${\mathfrak S}_w\mapsto{\mathfrak S}_{\overline{w}}$
on $H^*(\mbox{\em Flags}({\mathbb C}\,^n))$ interchanges
$e(\alpha)$ and $h(\alpha)$.

\section*{Acknowledgments}
We thank Mark Shimozono and Richard Stanley for helpful comments.
The second author is grateful to the  hospitality
of  Universit\'e Gen\`eve and INRIA Sophia-Antipolis,
where portions of this paper were developed and written.

\end{document}